\documentclass[11pt, reqno]{amsart}

\usepackage{amssymb,latexsym,amsmath,amsfonts,amsthm}
\usepackage{mathrsfs}
\usepackage{enumitem}
\usepackage{hyperref}
\usepackage{comment}

\usepackage[usenames]{color}
\definecolor{DPurple}{rgb}{0.46,0.2,0.69}

\allowdisplaybreaks

\numberwithin{equation}{section}

\theoremstyle{definition}
\newtheorem{definition}{Definition}[section]

\theoremstyle{remark}
\newtheorem{remark}[definition]{Remark}

\theoremstyle{plain}
\newtheorem{theorem}[definition]{Theorem}
\newtheorem{result}[definition]{Result}
\newtheorem{lemma}[definition]{Lemma}
\newtheorem{proposition}[definition]{Proposition}
\newtheorem{example}[definition]{Example}
\newtheorem{corollary}[definition]{Corollary}

\newcommand{\eps}{\varepsilon}

\newcommand{\ep}{\boldsymbol{\epsilon}}

\newcommand{\OM}{\Omega}
\newcommand{\bdy}{\partial}
\newcommand{\D}{\mathbb{D}}
\newcommand{\ome}{\omega}

\newcommand{\opn}{\mathscr{O}}
\newcommand{\opu}{\mathcal{U}}
\newcommand{\Util}{\widetilde{\mathcal{U}}}
\newcommand{\plDom}{\mathcal{D}(\beta,\eps)}

\newcommand{\smoo}{\mathcal{C}}

\newcommand{\Dini}{\mathcal{C}^{1,\,{\rm Dini}}}

\newcommand{\Ftil}{\widetilde{F}}
\newcommand{\Fhat}{\widehat{F}}

\newcommand{\bcdot}{\boldsymbol{\cdot}}

\newcommand{\btl}{\blacktriangleleft}
\newcommand{\nrml}{\eta_{\xi}}
\newcommand{\cluster}{C(F,p)}
\newcommand\domseq[1]{z_{\nu}^{\raisebox{-1pt}{$\scriptstyle{{#1}}$}}}
\newcommand\targseq[1]{w_{\nu}^{\raisebox{-1pt}{$\scriptstyle{{#1}}$}}}
\newcommand{\bdynbd}{\mathscr{V}}

\newcommand{\Z}{\mathbb{Z}}

\newcommand{\Cn}{\mathbb{C}^n}
\newcommand{\CC}{\mathbb{C}^2}
\newcommand{\C}{\mathbb{C}} 
\newcommand{\R}{\mathbb{R}}

\newcommand{\re}{{\sf Re}}
\newcommand{\im}{{\sf Im}}

\newcommand{\Sdist}{\delta_{D}}
\newcommand{\Tdist}{\delta_{\Omega}}
\newcommand{\TKobm}{k_{\Omega}}
\newcommand{\DKobd}{K_D}
\newcommand{\TKobd}{K_{\Omega}}
\newcommand{\DiscRadi}{\delta_{\Omega}(z;v)}

\begin{document}

\title[Local extension of proper holomorphic maps]{Local continuous extension of proper holomorphic \\ maps:
low-regularity and
infinite-type boundaries}

\author{Annapurna Banik}
\address{Department of Mathematics, Indian Institute of Science, Bangalore 560012, India}
\email{annapurnab@iisc.ac.in}

\begin{abstract}
  We prove a couple of results on local continuous extension of proper holomorphic maps $F:D \rightarrow
  \Omega$, $D, \Omega \varsubsetneq \Cn$, making local assumptions on $\bdy D$ and $\bdy \Omega$.
  The first result allows us to have much lower regularity, for the patches of $\bdy D, \bdy \Omega$
  that are relevant, than in earlier results. The second result (and a result closely related to
  it) is in the spirit of a result by
  Forstneri\v{c}--Rosay. However, our assumptions allow $\bdy \Omega$ to contain boundary points of
  infinite type. 
\end{abstract}


\maketitle

\vspace{-5.25mm}
\section{Introduction} \label{sec:intro}
In this paper, we present some results on local continuous extension of proper holomorphic
maps $F:D \rightarrow \Omega$, $D, \ \Omega \varsubsetneq  \Cn$, given local assumptions on $\bdy D$
and $\bdy \Omega$. Our results are motivated by the well-known work of Forstneri\v{c}--Rosay 
\cite{forstneric-rosay:1987}. There is also extensive literature on the problem of \emph{global}
extension of proper holomorphic maps (see \cite{forstneric:1993} and the references therein), 
from which too we take a few cues.
\smallskip

With $F, \ D$, and $\Omega$ as above, let $p \in \bdy D$ and let $\cluster$ be the
\emph{cluster set} of $F$ at $p$, defined as:
$$
  \cluster\,:=\, \big\{ w \in \Cn: \exists \text{ a sequence } \{z_{\nu}\}\subset D \ \text{such that} \, 
  \lim_{\nu \rightarrow \infty} z_{\nu} \,=\,p \ \text{and} \, \lim_{\nu \rightarrow \infty}
  F(z_{\nu})\,=\, w \big\}.
$$
If $\Omega$ is bounded, then $\cluster \neq \emptyset$. $F$ being proper, $\cluster
\subset \bdy \Omega$. 
In \cite{forstneric-rosay:1987}, the
``local assumptions'' alluded to above are imposed on $\bdy D \cap U$ and $\bdy \Omega \cap W$,
where $U$ and $W$ are neighbourhoods of $p$ and $q$ respectively, where
$q \in \cluster$. Since we will have occasion to mention the main result
in \cite{forstneric-rosay:1987} let us state it (also see \cite{berteloot:1992, sukhov:1994, mahajan:2012}
for results of a similar flavour as the following).

\begin{result}[{\cite[Theorem~1.1]{forstneric-rosay:1987}}] \label{r:F-R_main}
Let $D$ and $\Omega$ be domains in ${\C}^n$, $\Omega$ bounded, and let
$F: D \rightarrow \Omega$ be a proper holomorphic map. Let $p\in \bdy D$. Assume that there
is a continuous, negative plurisubharmonic function $\rho$ on $D$ and a neighbourhood $U$ of $p$ 
such that $\bdy D \cap U$ is a $\smoo^{1+\eps}$ submanifold of $U$ for some $\eps > 0$, and
$\rho(z) \geq -\delta_D (z)$ for all $z\in U \cap D$. If the cluster set $C(F,p)$ contains
a point $q$ at which $\bdy \OM$ is strongly pseudoconvex, then $F$ extends continuously to
$p$.
\end{result}
With $W$ and $q$ as above, Result~\ref{r:F-R_main} assumes $\bdy \Omega \cap W$ to be strongly pseudoconvex. Later results\,---\,\cite{berteloot:1992, sukhov:1994}, for instance\,---\,weaken
this requirement to admit certain families of domains $\Omega$ such that $\bdy \Omega$ is of finite type
at each point of $\bdy \Omega \cap W$. In this paper, we wish to extend this analysis and, in contrast,
focus on the situation when, among other things, $\bdy \Omega \cap W$ is not even $\smoo^1$-smooth
or, if it is smooth, then $q$ is of infinite type; see Examples~\ref{ex:Lipz_int-cone} and~\ref{ex:Dini_C2-cvx},
respectively.
\smallskip

Concerning the notation $\delta_D$ in Result~\ref{r:F-R_main}: given an open set
$D \varsubsetneq \C^n$ and $z\in D$, we write $\delta_D(z) := {\rm dist}(z, \C^n\setminus D)$.
The idea of the proof in \cite{forstneric-rosay:1987} (although not the details thereof)
goes back to Vormoor \cite{vormoor:1973}. The technique used in \cite{forstneric-rosay:1987} relies on
the classical Hopf Lemma for plurisubharmonic functions, which imposes constraints on the regularity
of $\bdy \Omega \cap W$. For $F, \ D, \ \Omega$, and $p$ as above,
and given a sequence $\{z_{\nu}\} \subset D$ such that $z_{\nu} \to p$, some form of a Hopf-type lemma
is the simplest tool to control $\{F'(z_{\nu})\}$\,---\,provided $\bdy \Omega$ is at least
$\smoo^2$ near $\{F(z_{\nu})\}$ (see \cite{henkin:1973, pinchuk:1974, diederich-fornaess:1979}, for instance). 
We explore a different paradigm from the one in \cite{forstneric-rosay:1987} that allows us to
greatly lower the regularity of $\bdy D \cap U$ and $\bdy \Omega \cap W$. Loosely speaking,
instead of picking a point in $\cluster$ and imposing conditions on $\bdy \Omega$ near it,
we consider an \emph{interior} condition (i.e., on a suitable open set in $\Omega$) which is
notably less restrictive; see Theorem~\ref{th:non_smooth_bdy_extn}. Before we present this
result, a few definitions:

\begin{definition}
A function $\ome: ([0,\infty),0) \rightarrow ([0,\infty),0) $ is said to satisfy a 
\emph{Dini condition} if $\ome$ is monotone increasing and $\int_{0}^{\eps}
r^{-1}\ome(r) dr <\infty$ for some (hence for any) $\eps>0$.
\end{definition}\label{defn:Dini}

\begin{definition} \label{defn:Lipz}
Let $D \subset \Cn $ be a domain, $p \in \bdy D $, and $U$ be a neighbourhood of $p$.
Let $S \,:=\, \bdy D \cap U$.
\begin{enumerate}[leftmargin=25pt]
  \item We say that $S$ is a \emph{Lipschitz submanifold of $U$} if for each $\xi \in S$,
  there exists a neighbourhood $\mathcal{N}_\xi$ of $\xi$, $\mathcal{N}_\xi \subset U$, a
  unitary map ${\sf U}_{\xi}$, a constant $r_\xi>0$, and a Lipschitz function
  $\varphi_{\xi}: B^{n-1}(0,r_\xi) \times (-r_\xi, r_\xi) \rightarrow \R$ such that, writing
  the affine map $z \mapsto {\sf U}_{\xi}(z- \xi)$ as $\sf{U}^{\xi}$, we have
  \begin{align}
    {\sf U}^{\xi}(\mathcal{N}_\xi \cap D)  \subset\ \{ &(Z',Z_n)\in {\C}^{n-1} \times \C:
    \im(Z_n) > \varphi_\xi (Z', \re(Z_n)), \ \|Z'\|< r_\xi \notag \\ 
    &\,\text{and} \
    | \re (Z_n) | < r_\xi \}, \label{eqn:dom} \\
    {\sf U}^{\xi}(\mathcal{N}_\xi \cap \bdy D) =\ \{ &(Z',Z_n)\in {\C}^{n-1} \times \C:
    \im(Z_n) = \varphi_\xi (Z', \re(Z_n)), \ \|Z'\|< r_\xi \notag \\
    &\,\text{and} \ | \re (Z_n) | < r_\xi \}. \label{eqn:bdy}
  \end{align}
  \item (Cf.~Kukavica--Nystr{\"o}m \cite{kukavica-nystrom:1998}, for instance) We say
  that $S$ is a \emph{$\Dini$ submanifold of $U$} if
  there exists a function $\omega: ([0, \infty), 0)\to
  ([0, \infty),0)$ satisfying a Dini condition, and for each $\xi \in S$,
  there exist $\mathcal{N}_\xi$, ${\sf U}^{\xi}$, and $r_\xi>0$ as described by (1) above,
  and a $\smoo^1$ function $\varphi_{\xi}: B^{n-1}(0,r_\xi) \times (-r_\xi, r_\xi)
  \rightarrow \R$ such that
  $$
    \|\nabla \varphi_{\xi}(x) - \nabla \varphi_{\xi}(y)\|
    \,\leq\,\omega( \|x - y\| ) \quad \forall x, y \in B^{n-1}(0,r_\xi) \times (-r_\xi, r_\xi),
  $$
  and ${\sf U}^{\xi}(\mathcal{N}_\xi \cap D)$ and ${\sf U}^{\xi}(\mathcal{N}_\xi \cap \bdy D)$
  are described by \eqref{eqn:dom} and \eqref{eqn:bdy}, respectively.
\end{enumerate}
When $n=1$, the expressions $B^{n-1}(0,r_\xi) \times (-r_\xi, r_\xi)$ and ${\C}^{n-1} \times \C$
must be read as $(-r_\xi, r_\xi)$ and $\C$, respectively. The conditions \eqref{eqn:dom} and
\eqref{eqn:bdy} must then be read \emph{mutatis mutandis}.
\end{definition}

Before we give the next definition, let us fix some notations. An \emph{open right circular
cone with aperture $\theta$} is the following open set 
$$
  \Gamma(v,\theta) \,:=\, \{ z \in {\C}^n : \re \, \langle z,v \rangle > \cos (\theta /2)\,  \|z\| \},
$$
where $v \in {\C}^n$ is a unit vector and $\theta \in (0, \pi)$ 
(here $\langle \bcdot\,, \bcdot \rangle$ denotes the standard Hermitian inner product
on ${\C}^n$). Given $z_0 \in \Cn$, the \emph{axis} of the (translated) cone 
$ z_0 + \Gamma (v, \theta) $ is the ray $\{z_0 + tv: t>0\}$.

\begin{definition} \label{defn:int_cone}
Let $D \subset \Cn$ be a domain and let $W$ be an open set such that $\bdy D \cap W$ is a
non-empty $\bdy D$-open set. We say that $D$ satisfies a \emph{uniform interior cone
condition in $W$} if, given any open set $\opu \subset D$ such that $\opu \Subset
W$, there exist constants $r_{\opu}>0$ and $\theta_{\opu} \in (0,\pi)$ such that for each
$w \in \opu$, there exists $\xi_{w} \in \bdy D \cap W$ and a unit vector $v_w$ such that
\begin{itemize}
  \item $w$ lies on the axis of $\xi_w + \Gamma(v_w, \theta_{\opu})$.
  \smallskip
  \item $w \in (\xi_w + \Gamma(v_w, \theta_{\opu})) \cap B^n(\xi_w, r_{\opu}) \subset W \cap D$.
\end{itemize}
\end{definition}

Given a domain $D \subset \C^n$, $k_{D} : D\times \C^n \rightarrow [0, \infty)$
will denote the Kobayashi pseudometric for $D$.
We are now in a position to state our first theorem. 

\begin{theorem} \label{th:non_smooth_bdy_extn}
Let $D$ and $\Omega$ be domains in ${\C}^n$ and let $F: D \rightarrow \Omega$
be a proper holomorphic map. Let $p\in \bdy D$. Assume that there is a continuous,
negative plurisubharmonic function $\rho$ on $D$, a neighbourhood $U$ of $p$, and a constant
$s\in (0,1]$ such that $\bdy D \cap U$ is a Lipschitz submanifold of $U$, and
$\rho(z) \geq -(\delta_D(z))^s$ for all $z\in U \cap D$. Suppose there exists a
neighbourhood $U^*$ of $p$, $U^* \Subset U$, and an open set $W$ such that 
$\bdy \Omega \cap W \neq \emptyset$, and such that
\begin{itemize}
    \item $F(U^* \cap D)\Subset W$, and
    \smallskip
    \item $\Omega$ satisfies a uniform interior cone condition in $W$.
\end{itemize}
Suppose there exists a function $M: ([0, \infty), 0)\to
([0, \infty),0)$ satisfying a Dini condition so that
\begin{equation}\label{eqn:Kobm_bound}
  k_{\Omega}(w; v) \geq \|v\|/M(\delta_{\Omega}(w)) \quad \forall (w, v) 
  \in (W \cap \Omega)\times {\C}^n.
\end{equation}
Then, there exists a $\bdy D$-neighbourhood $\bdynbd$ of $p$ such that $F$ extends
to $D \cup \bdynbd$ as a continuous map.
\end{theorem}

One may ask whether, given that $\Omega$ in Theorem~\ref{th:non_smooth_bdy_extn}
is assumed to satisfy a uniform interior cone condition in $W$, one also requires the
condition \eqref{eqn:Kobm_bound}. There is a vital point related to this, which
is best discussed after we state our next theorem and prove Theorem~\ref{th:non_smooth_bdy_extn};
see Remark~\ref{rem:Kobm}.
\smallskip


The proof of our next result follows many of the techniques used in
\cite{forstneric-rosay:1987}. Even so, with $D,\ p$, and $U$ as above, one is able to admit
$\bdy D \cap U$ having lower regularity than in \cite{forstneric-rosay:1987}. The
argument used in \cite{forstneric-rosay:1987} is \emph{very delicate}, and is necessitated
by an aspect of the hypothesis of Forstneri{\v{c}}--Rosay: i.e., picking a point
in $\cluster$ and imposing conditions on $\bdy \Omega$ near it (compare the hypotheses of
Result~\ref{r:F-R_main} and Theorem~\ref{th:non_smooth_bdy_extn}). Given such a hypothesis,
the proof relies on certain intrinsic constants matching up precisely. This is why, as hinted
at earlier, one needs to assume $\mathcal{C}^2$ regularity near a point $q \in \cluster$: which
is done both in Result~\ref{r:F-R_main} and in Theorem~\ref{th:Dini_C2-cvx}.
Our greatest departure from \cite{forstneric-rosay:1987} involves a concept\,---\,namely,
local log-type convexity\,---\,introduced by Liu--Wang in
\cite{liu-wang:2021}. The assumption of $\Omega$ being log-type convex near $q$ (which
substitutes the assumption of strong pseudoconvexity in Result~\ref{r:F-R_main})
is broad enough to admit, on the one hand, domains $\Omega$ such that
$\bdy \Omega \cap W$ is a $\smoo^2$ submanifold and, on the other hand, $\Omega$
such that $q$ is of infinite type. Quantitatively,
this condition admits useful lower bounds for the Kobayashi distance for $\Omega$,
which is the delicate part of proving Theorem~\ref{th:Dini_C2-cvx}.
\smallskip

Before we present Theorem~\ref{th:Dini_C2-cvx}, we need to define log-type convexity.
But first, we must introduce a
quantity that is very natural in the context of bounded convex domains. For instance,
by a theorem of Graham \cite{graham:1990, graham:1991} (see Result~\ref{r:Graham_cvx}
below), this quantity provides an optimum estimate for the Kobayashi metric for such a domain. Let
$D$ be a bounded convex domain in $\Cn$. For each $z \in D$ and $v \in \Cn \setminus
\{0\}$,~define 
$$
  \Sdist(z;v) \,:=\, \sup \bigg\{r>0: \bigg(z+ (r\D) \frac{v}{\|v\|} \bigg) 
  \subset D\bigg\}.
$$

\begin{definition}[Liu--Wang, {\cite[Definition~1.1]{liu-wang:2021}}]  \label{defn:L-T-C}
A bounded convex domain $D \subset \Cn, \ n \geq 2$, is called \emph{log-type convex} if
there are constants $C, \ \nu >0$ such that
\begin{equation} \label{eqn:L-T-C}
  \Sdist(z;v) \leq \frac{C}{{\big|\log \Sdist(z) \big|}^{1 + \nu}} \quad \forall z \in D \ \text{and}
  \ \forall v \in \Cn \setminus \{0\} .
\end{equation}
\end{definition}

\begin{theorem}  \label{th:Dini_C2-cvx}
Let $D$ and $\Omega$ be domains in ${\C}^n, \ n\geq 2, \ \Omega$ bounded, and 
let $F: D \rightarrow \Omega$
be a proper holomorphic map. Let $p\in \bdy D$ and $q \in C(F,p)$. Assume that there is a
continuous, negative plurisubharmonic function $\rho$ on $D$ and a neighbourhood $U$ of $p$ 
such that $\bdy D \cap U$ is a $\Dini$ submanifold of $U$, and $\rho(z) 
\geq - \delta_D (z)$ for all $z\in U \cap D$. Suppose there exists a neighbourhood $\opn$
of $q$ such that 
\begin{itemize}
  \item $\bdy \Omega \cap \opn$ is a $\mathcal{C}^2$ submanifold of $\opn$, and
  \smallskip
  \item $\opn \cap \Omega$ is log-type convex.
\end{itemize}
Then, there exists a $\bdy D$-neighbourhood $\bdynbd$ of $p$ such that $F$ extends
to $D \cup \bdynbd$ as a continuous map.
\end{theorem}

\begin{remark}
The reason Theorem~\ref{th:Dini_C2-cvx} has been stated for domains in ${\C}^n$ with $n \geq 2$
is because its proof relies crucially on results about log-type convex domains (see
Section~\ref{sec:lemmas_Dini_C2-cvx}), which are defined as domains in ${\C^n}$, $n \geq 2$.
Note, though, that \emph{any} bounded convex domain in $\C$ satisfies the defining inequality
\eqref{eqn:L-T-C} for a log-type convex domain. However, the proofs of the above-mentioned results
have been written with domains in higher dimensions in mind\,---\,additional arguments would be
needed in the planar case. For this reason, we restrict Theorem~\ref{th:Dini_C2-cvx} to
$n \geq 2$. Furthermore, one suspects that Theorem~\ref{th:Dini_C2-cvx}, suitably restated,
is already known for planar domains.
\end{remark}

We must mention a connection between Theorem~\ref{th:non_smooth_bdy_extn} and the proof
of Theorem~\ref{th:Dini_C2-cvx}. Note that the conclusion of Theorem~\ref{th:Dini_C2-cvx} is stronger
than that of Result~\ref{r:F-R_main}. Now, the former conclusion can, in principle, be deduced under
the assumptions in Result~\ref{r:F-R_main} too. This requires an auxiliary argument alluded to in
\cite{forstneric-rosay:1987} (see page 241)\,---\,by which one may deduce H{\"o}lder
continuity of the extension with exponent $1/2$. 
That specific argument does not, in general, work in the
context of Theorem~\ref{th:Dini_C2-cvx}. Instead, it turns out that once we get a continuous extension
of $F$ to $p$ (for Theorem~\ref{th:Dini_C2-cvx}), we are able to use
Theorem~\ref{th:non_smooth_bdy_extn} and get the stronger conclusion. We can also estimate
the modulus of continuity of the extension given by Theorems~\ref{th:non_smooth_bdy_extn}
and~\ref{th:Dini_C2-cvx}. As both theorems admit domains $\Omega$ such that $\bdy \Omega \cap W$
contains points of infinite type (see Remark~\ref{rem:Kobm} and Example~\ref{ex:Dini_C2-cvx}), one
\textbf{cannot} expect these extensions, in general, to be H{\"o}lder. This leads to the much more
technical discussion of the modulus of continuity of the extension, which we omit.
\smallskip

We should also emphasise that\,---\,with $D$, $U$, and $p\in \bdy{D}$ as in all the results above\,---\,the observations just made hold true for $\bdy{D}\cap U$ having much lower regularity than in any
earlier result on the present theme.
\smallskip

Readers familiar with the essence of the argument of Forstneri\v{c}--Rosay in
\cite{forstneric-rosay:1987} and with the results in \cite{bracci-nikolov-thomas:2022}
by Bracci \emph{et al.} might ask whether the conclusions of Theorem~\ref{th:Dini_C2-cvx}
may be obtainable under weaker conditions on $\opn \cap \Omega$. We shall address this in
Section~\ref{sec:second_th_proof}: see Remark~\ref{rem:cvx_visi} and Theorem~\ref{th:Dini_C2-cvx-visi}.%
\smallskip

Given any two domains $D,\,\Omega \subseteq \C^n$, $n \geq 2$, it is rare for the pair $(D, \Omega)$ to
admit a proper holomorphic map $F : D \rightarrow \Omega$. Since the domains $D,\,\Omega$ in
Theorems~\ref{th:non_smooth_bdy_extn} and~\ref{th:Dini_C2-cvx} must satisfy several conditions,
the question arises: are there any domains $D,\,\Omega \subseteq \C^n$ that satisfy these conditions
\textbf{and} admit a non-trivial proper holomorphic map $F : D \rightarrow \Omega$ when
$n \geq 2$? We provide examples
of such domains in Section~\ref{sec:ex}. Recall: the conditions in Theorem~\ref{th:Dini_C2-cvx}
allow $\bdy \Omega$ to be of infinite type at $C(F, p)$. Section~\ref{sec:ex} provides an
example of $(D, \Omega)$ where $\Omega$ has the latter property and there exists a non-trivial proper
holomorphic $D \rightarrow \Omega$ map. As for the proofs of Theorems~\ref{th:non_smooth_bdy_extn}
and~\ref{th:Dini_C2-cvx}: they are presented in Sections~\ref{sec:first_th_proof}
and~\ref{sec:second_th_proof}, respectively. 
\medskip

\section{Two examples}  \label{sec:ex}
In this section, we discuss the examples mentioned above. But
we first explain the notation used below and in later sections (some of which has also
been used without clarification in Section~\ref{sec:intro}).
\smallskip

\subsection{Common notations}
\begin{enumerate}[leftmargin=25pt]

  \item For $v \in {\R}^N $, $\|v\|$ denotes the
  Euclidean norm. For any $a \in {\R}^N$ and $B \subset {\R}^N$, we write $\text{dist}(a, B)\,:=\,
  \inf \{\|a-b\|: b \in B\}$. 
  \smallskip
  \item Given a point $z \in \Cn$ and $r>0$, $B^n(z,r)$ denotes the open Eulidean ball with radius
  $r$ and center $z$. For simplicity, we write $\D \,:=\, B^1(0,1)$.
  \smallskip
  \item Given a domain $D \subset \Cn$, $K_D:D \times D \rightarrow [0, \infty)$ denotes the
  Kobayashi distance for $D$.
\end{enumerate}

\subsection{Examples} We are now in a position to present the examples referred to several times
in Section~\ref{sec:intro}. To this end, we need the
following result by Sibony:

\begin{result}[paraphrasing {\cite[Proposition~6]{sibony:1981}}]\label{r:Kob_low_bd_Sibony}
Let $D \subset \Cn$ be a domain and $z \in D$. There exists a uniform constant
$\alpha>0$ (i.e., it does not depend on $z$) such that if $u$ is a negative plurisubharmonic
function on $D$ that is of class $\mathcal{C}^2$ in a neighbourhood of $z$ and satisfies 
$$
  \langle v, (\mathfrak{H}_{\C}u)(z) v \rangle \geq c \|v\|^2 \quad \forall v \in \Cn,
$$
where $c$ is some positive constant, then 
$$
  k_D(z;v) \geq \Big( \frac{c}{\alpha} \Big)^{1/2}\frac{\|v\|}{|u(z)|^{1/2}}.
$$
\end{result}

\noindent{Here $\langle \bcdot \,, \bcdot \rangle$ denotes the standard Hermitian inner product and}
$\mathfrak{H}_{\C}$ denotes the complex Hessian.

\begin{example} \label{ex:Lipz_int-cone}
An example demonstrating that there exist domains $D,\ \Omega \subset \CC,$ and a proper holomorphic map
$F:D \rightarrow \Omega$ such that $D, \ \Omega,$ and $F$ satisfy the conditions stated in Theorem~\ref{th:non_smooth_bdy_extn}.
\end{example}

\noindent{Let us define 
\begin{align*}
  D\,&:=\, \{(z,w) \in \CC: |z|^2 + |w| < 1 \},\\
  \Omega\,&:=\, \{(z,w)\in \CC: |z|+ |w| < 1 \},
\end{align*}
and $F(z,w) := (z^2,w)$ for $(z,w) \in D.$
Then $F:D \rightarrow \Omega$ and it is clear that $F$
is proper. Finally, let $p=(1,0) \in \bdy D$. That the stated conditions are satisfied will
be discussed in two steps.} 

\medskip
\noindent{{\textbf{Step 1.}} We shall show that there exists a constant $\widetilde{C}>0$
such that, writing 
$$
  \rho(z,w) \,:=\, \widetilde{C}\, \big(|z|^2 + |w| -1 \big), \ (z,w) \in D,
$$
$\rho$ satisfies the desired estimate in $U \cap D$, where $U\,:=\, 
\{(z,w)\in \CC:9/10<|z|<11/10,~|w|<~1\}$. Since $\bdy D$ is \textbf{not} smooth around $p$,
this task will be slightly involved. Note that} $\rho$ is a continuous, negative, plurisubharmonic
function on $D$. 
\smallskip

Since $D$ is Reinhardt, it is elementary to show that if $(z,w) \in D, \ \theta_z, \ \theta_w \in
\R$ are such that $z=|z|e^ {i\theta_z},\ w= |w|e^{i \theta_w}$, and $(\zeta, \eta) \in \bdy D$ are 
such that $\|(z,w)-(\zeta,\eta) \| \,=\, \Sdist(z,w)$, then
\begin{itemize}
  \item $\exists (\zeta, \eta) \in \bdy D$ with the above property such that $\zeta=|\zeta|e^{i\theta_z}$
  and $\eta = |\eta| e^{i\theta_w}$.
  \smallskip
  
  \item For $(\zeta, \eta)$ as described by the above bullet point, $\big\| \big(|z|,|w| \big)-
  \big(|\zeta|, |\eta| \big)\big\|=\Sdist \big(|z|,|w| \big).$
\end{itemize}
Due to the above, it suffices to show that 
\begin{align}   \label{eqn:first_quad_Sbound}
  \rho \big(|z|,|w| \big) \geq - \Sdist \big(|z|,|w| \big) \quad \forall (z,w): |z|^2 + |w|<1, \ 
  9/10<|z|<1, \ |w|< 1/10.
\end{align}

Clearly, to establish \eqref{eqn:first_quad_Sbound}, we need to estimate the infimum of the set
$\mathcal{S}(x_0,y_0)\,:=\, \{\|(x,y)-(x_0,y_0)\|: x^2 + |y| = 1 \}$, having fixed $(x_0,y_0)$ such 
that $9/10< x_0 <1, \ 0 \leq y_0 < 1/10$, and ${x_0}^2 + y_0 <1$. It is elementary to show that, as 
the curve $x^2 + |y| =1$ has a non-smooth point at $(1,0)$ and bounds a convex region,
$\mathcal{S}(x_0,y_0)$ must attain its minimum in the set $ \{ (x,y)\in \R \times (0, \infty):
x^2+y=1 \}$, which is a smooth curve $\mathscr{C}$. So, we can apply the method of Lagrange multipliers to deduce
that if $\mathcal{S}(x_0,y_0)$ attains its minimum at $(X,Y)$, then $(X,Y)$ satisfies
\begin{align}
  2X^3+(2y_0 -1) X - x_0 &=0,\label{eqn:L_mltp_1} \\
  \frac{X-x_0}{Y-y_0} &= 2X. \label{eqn:L_mltp_2}
\end{align}

To estimate min $\mathcal{S}(x_0,y_0)$, we define an auxiliary function $\phi: \R \rightarrow \R, \ 
\phi(x) \,:=\, 2x^3+ (2y_0-1)x -x_0$. Independent of the choice of $y_0$, $\phi$ is strictly increasing
on $[1/ \sqrt{6}, \infty)$. By the nature of the curve $\mathscr{C}$ and the choice of $x_0$, it follows
that $9/10<x_0<X<1$. Thus, in order to locate $X$, it suffices to examine $\phi|_{\scriptscriptstyle{[9/10,1]}}$. 
Let $C\,:=\, \sup \{6x^2+(2y-1): x \in [9/10,1], \ y \in [0,1/10] \}$. Since $\sup_{\scriptscriptstyle{x\in [9/10,1]}}
\phi'(x)\leq C$, writing $A=\phi(9/10), \ B= \phi(1)$ we have
\begin{align}   \label{eqn:auxl_eqn}
  |\phi(x)-\phi(x')| &\leq C|x-x'| \quad \forall x,\ x' \in [9/10,1] \notag \\ 
  \Rightarrow \big|\phi^{-1}(a)-\phi^{-1}(a') \big| &\geq C^{-1}|a-a'| \quad \forall a,\ a' \in [A,B]. 
\end{align}
Thus, by \eqref{eqn:L_mltp_1} and \eqref{eqn:auxl_eqn}, since $X$ is a root of \eqref{eqn:L_mltp_1},
$$
  |X-x_0|=\big|\phi^{-1}(0)-\phi^{-1}(\phi(x_0)) \big| \geq C^{-1}|\phi(x_0)|=
  \frac{2x_0}{C}\,\big|x_0^2+y_0-1 \big|.
$$
From the above inequality, we deduce 
\begin{align}   \label{eqn:psh_bd}
  \min \mathcal{S}(x_0,y_0)=\sqrt{(X-x_0)^2+(Y-y_0)^2} \geq |X-x_0| \geq \frac{9}{5C} \big|x_0^2+y_0-1 \big|.
\end{align}
We set $\widetilde{C}=\frac{9}{5C}$. Then, \eqref{eqn:psh_bd} implies $\rho (x_0,y_0) \geq 
- \min \mathcal{S}(x_0,y_0)$. Recalling the purpose of the set $\mathcal{S}(x_0,y_0), \ \min \mathcal{S}(|z|,|w|)
=\Sdist(|z|, |w|)$ for the above mentioned constraints. Thus, the last inequality gives us
\eqref{eqn:first_quad_Sbound}. By the discussion preceding \eqref{eqn:first_quad_Sbound},
$$
  \rho(z,w) \geq - \Sdist(z,w) \quad \forall (z,w) \in U \cap D.
$$

\noindent{\textbf{Step 2.} We shall now show that $\TKobm$ satisfies the desired estimate in $(W \cap \Omega)
\times \CC$ where $W\,:=\, B^2((1,0),\delta),\ \delta>0$ is small enough\,---\,with  
$M(t)\,:=\,\Tilde{c}\,t^{1/2},\ t \geq 0$, for some $\Tilde{c}>0$}.
\smallskip

Clearly, $u(z,w)\,:=\, |z| + |w| -1, \ (z,w) \in \Omega $, is negative, plurisubharmonic on $\Omega$.
Note that $u$ is smooth at a given point $(z,w) \in W\cap \Omega$ if and only if $w \neq 0$. To 
establish the required lower bound on $\TKobm$ at smooth points of $u$, we will use
Result~\ref{r:Kob_low_bd_Sibony}. For the non-smooth points, we will use the convexity 
of $\Omega$ and apply the estimate by Graham as given by Result~\ref{r:Graham_cvx}.
\smallskip

Since $\Omega$ is Reinhardt, it satisfies the same properties as described by the bullet points (prior to
\eqref{eqn:first_quad_Sbound}) in Step~1. Thus, using coordinate geometry, it is elementary to see that 
\begin{equation}    \label{eqn:Tdist}
  \Tdist(z,w)\,=\,\frac{1-|z|-|w|}{\sqrt{2}}\,=\,\frac{|u(z,w)|}{\sqrt{2}}.
\end{equation}

Let $(z,w) \in W \cap \Omega$ with $w \neq 0$. It is easy to compute that
$$
  \langle v, (\mathfrak{H}_{\C}u)(z,w)v \rangle = \frac{1}{4} \bigg( \frac{|v_1|^2}{|z|} + 
  \frac{|v_2|^2}{|w|}\bigg) \geq \frac{\|v\|^2}{4} \quad \forall v=(v_1,v_2) \in \CC.
$$
Then, by Result~\ref{r:Kob_low_bd_Sibony} and \eqref{eqn:Tdist}, there exists $\beta>0$ such that 
\begin{equation}    \label{eqn:smooth_pts}
  \TKobm((z,w);v) \geq \frac{\beta \|v\|}{(\Tdist(z,w))^{1/2}} \quad \forall (z,w)\in W \cap \Omega \
  \text{with} \ w \neq 0 \ \text{and} \ \forall v \in \CC.
\end{equation}

Now, let $(z,w) \in W \cap \Omega$ with $w\,=\,0$ and let $\ v=(v_1,v_2) \in \CC$. Write $z\,=\,|z| e^{i\theta_z},\ 
v_k\,=\,|v_k|e^{i\theta_k}$ where $\theta_z,\ \theta_k \in \R,\ k\,=\, 1,\ 2$. By the invariance of $\TKobm$ under
the automorphism $\Omega~\ni~(z,w) \mapsto~(e^{-i\theta_z}z, e^{-i(-\theta_z + \theta_1-\theta_2)}w \big)$, we get
\begin{align} \label{eqn:aut_Graham}
  \TKobm((z,0);v)\,=\,\TKobm 
  \big( (|z|,0);e^{i(\theta_1-\theta_z)}\big(|v_1|,|v_2|\big) \big)
  &\geq \frac{\|v\|}{2 \, \Tdist \big( (|z|,0);e^{i(\theta_1-\theta_z)}\big(|v_1|,|v_2| \big) \big)}
  \notag \\
  &=\,\frac{\|v\|}{2 \, \Tdist \big( (|z|,0);\big(|v_1|,|v_2| \big) \big)},
\end{align}
where the inequality is due to Result~\ref{r:Graham_cvx}. It is easy to see that 
$$
  \Tdist \big( (|z|,0);\big(|v_1|,|v_2| \big) \big) \leq \big\| (|z|,0)-(1,0) \big\| 
  \,=\, \sqrt{2} \, \Tdist(z,0).
$$
Hence, \eqref{eqn:aut_Graham} implies 
\begin{equation}    \label{eqn:non_smooth_pts}
  \TKobm ((z,0);v) \geq \frac{\|v\|}{2\sqrt{2} \,\Tdist(z,0)} \quad \forall z:(z,0) \in W\cap \Omega \
  \text{and} \ \forall v \in \CC.
\end{equation}
If we set $\Tilde{c}\,=\,$max$\big( 1/ \beta, 2\sqrt{2} \big) $, then from \eqref{eqn:smooth_pts} and
\eqref{eqn:non_smooth_pts} we get the estimate
$$
  \TKobm((z,w);v) \geq \|v\|/M(\Tdist(z,w)) \quad \forall (z,w) \in W\cap \Omega \ \text{and} \ \forall v \in \CC
$$
for the $M$ introduced above.
\smallskip

All that remains is to produce a neighbourhood $U^*$ of $(1,0)$ such that $(U^* \cap D)$ is
mapped by $F$ as desired. Now, the main point of this example is\,---\,given that proper holomorphic
maps between a given pair of domains in $\C^n$ are rare when $n \geq 2$\,---\,to show that there exist domains $D$
and $\Omega$ satisfying the respective \emph{geometric} conditions imposed by Theorem~\ref{th:non_smooth_bdy_extn}
that admit a proper non-trivial holomorphic map $F : D \rightarrow \Omega$. In the present example, $F$ is holomorphic
on $\C^2$. Clearly, it extends continuously to a $\bdy D$-neighbourhood, say $\bdynbd$, of
$(1,0)$. The existence of the desired $U^*$ easily follows from the continuity of $F$ at
$p\,=\,(1,0)$. \hfill $\btl$
\smallskip

Recall the discussion in Section~\ref{sec:intro} on how the assumptions stated in
Theorem~\ref{th:Dini_C2-cvx} admit continuous extension of $F : D \rightarrow \Omega$ even if 
$q\in \bdy \Omega$\,---\,$D$, $\Omega$, $F$, and $q$ as in the statement of
Theorem~\ref{th:Dini_C2-cvx}\,---\,is a point of infinite type. This is very different
from the assumptions encountered in previous results on local extension of proper holomorphic
maps. For this reason, an example, showing that the result implicit in the first sentence of this paragraph
isn't vacuously true, is desirable. In discussing such an example, the following well-known result
(see, for instance, \cite[Chapter~3]{krantz:2001}) will be useful.

\begin{result}  \label{r:bdy_dist_defng_fn}
If $D$ is a domain in $\Cn, \ p \in \bdy D$, and $\bdy D$ is a $\mathcal{C}^2$ submanifold of some
neighbourhood $\mathcal{N}$ of $p$, then there exists an $\eps > 0$ such that $B^n(p, \eps) \Subset
\mathcal{N}$ and: 

\begin{itemize}
  \item $\Sdist (z) \,=\, \delta_{\mathcal{N} \cap D}(z) $ for all $z \in B^n(p, \eps) \cap D $.
  \smallskip
  
  \item If we define 
  $$
   \widetilde{\rho}(z) \,:=\,\begin{cases}
  					         -\Sdist(z), &\text{if $z \in \overline{D} $}, \\
  					         \delta_{(\Cn\setminus \overline{D})}(z), &\text{if $z \in \Cn \setminus \overline{D}$},
  					         \end{cases}
  $$
\end{itemize}
then $\widetilde{\rho} \in \mathcal{C}^2 \big( B^n (p,\eps)\big)$. Moreover, $\widetilde{\rho}$ is
a local defining function for $D$ on $B^n(p,\eps)$.
\end{result}

\begin{example} \label{ex:Dini_C2-cvx}
An example demonstrating that there exist domains $D,\ \Omega \subset \CC,$ and a proper holomorphic
map $F:D \rightarrow \Omega$ such that $D, \ \Omega,$ and $F$ satisfy the conditions stated in
Theorem~\ref{th:Dini_C2-cvx} and such that, with $p\in \bdy D$ and for $q \in C(F,p)$, $\bdy \Omega$
is of infinite type at $q$.
\end{example}

\noindent{Let us define the function $\varphi : [0, \infty) \rightarrow \R$ as
$$
 \varphi(x) \,:=\, \begin{cases}
  					\exp \big(\!\!-\frac{1}{\sqrt{x} }\big), &\text{if $x \neq 0$}, \\
  					
  					0, &\text{otherwise}.
  					\end{cases}
$$
Now, set
\begin{align*}
  D\,&:=\, \{(z,w) \in \CC: \re z > \varphi(|w|^2)\} \cap 
  \{(z,w) \in \CC:|z|^2 + |w|^4 <1 \},\\
  \Omega\,&:=\, \{(z,w) \in \CC: \re z > \varphi(|w|)\} \cap 
  B^2((0,0),1),
\end{align*}
and $F(z,w) := (z, w^2)$ for $(z,w) \in D$. Then $F:D \rightarrow \Omega$
and it is clear that $F$ is proper. Let 
$p \,=\, (0,0) \in \bdy D$ and $q \,=\, (0,0) \in C(F,p)$. Take $\opn \,:=\, B^2(q,r)$, 
where $r \in (0,1)$ is such that $\opn \cap \Omega$ is log-type convex (see \cite[Example~1.3]
{liu-wang:2021}). Consider the function $\rho : D \rightarrow \R$ defined as
$$ 
  \rho (z,w) \,:=\, \varphi(|w|^2) - \re z, \ (z,w) \in D.
$$
By definition, $\rho$ is smooth and negative in $D$. A computation gives
$$
  \frac{{\bdy}^2 \rho}{\bdy w \bdy \overline{w}}(z,w) \,=\,
  \begin{cases}
    {4}^{-1}  {|w|}^{-3} \exp \big(\!\!-\frac{1}{|w| }\big) \Big(\frac{1}{|w| } -1 \Big), &\text{if $w \neq 0$}, \\
    0, &\text{otherwise}.
  \end{cases}
$$
This implies that $\rho$ is a plurisubharmonic function on $D$.
\smallskip

We will now show that there exist a constant $C>0$ and a neighbourhood $U$ of $p$ such that the 
plurisubharmonic function $D \ni (z,w) \mapsto C \rho(z,w)$ satisfies the desired estimate 
in $U \cap D$.} Since $\bdy D$ is smooth near $p$, there exist a neighbourhood $\mathcal{N}$ of $p$ and
$\eps>0$ with $B^2(p, \eps) \Subset \mathcal{N} \Subset \{(z,w) \in \CC:|z|^2 + |w|^4 <1 \} $
such that the conclusion of Result~\ref{r:bdy_dist_defng_fn} holds and such that
$\rho|_{\mathcal{N}}$ is a local defining function for $D$ at $p \in \bdy D$. Now, since the function
$\widetilde{\rho}$ as in Result~\ref{r:bdy_dist_defng_fn} is a local defining function of $D$ on
$B^2(p, \eps)$, there exists a neighbourhood $\mathcal{N}'$ of $\bdy D \cap
B^2(p,\eps)$ 
and a positive function $H \in \mathcal{C}^{2}(\mathcal{N}')$ 
satisfying  
$$
  \widetilde{\rho}(z,w)\,=\, \rho(z,w) H(z,w) \quad \forall (z,w) \in \mathcal{N}'.
$$
So, fixing a neighbourhood $U$ of $p$ such that $U \Subset \mathcal{N}'$, write $C\,:=\,
\inf_{U}H >0$. Then,
$$
  C \rho(z,w) \geq - \Sdist(z,w) \quad \forall (z,w) \in U \cap D. 
$$
Finally, note that $\bdy D$ is of infinite type at $q=(0,0)$.
Hence, the triple $(D, \Omega, F)$ has the desired properties.  \hfill $\btl$
\medskip

\section{Preliminary analytic lemmas}   \label{sec:analytic_lemmas}
This short section is devoted to a few facts from analysis that would be needed
in the proofs of the theorems stated in Section~\ref{sec:intro}.

\begin{lemma} \label{l:Lipz_vertical}
Let $p\in \bdy D$ and $U, \ U^{*}$ be as in Theorem~ \ref{th:non_smooth_bdy_extn}. Then
(in the notation of Definition~\ref{defn:Lipz})
there exist a neighbourhood $V$ of $p,\ V \Subset \mathcal{N}_{p}  \cap U^{*}$ and a
constant $C>1$ such that
$$
  \Sdist(z) \leq Y({\sf U}^p(z)) \leq C \Sdist(z) \quad \forall z \in V\cap D
$$
where we define $ Y(Z',Z_n) \,:=\, \im (Z_n)-\varphi_p (Z', \re(Z_n))$ 
for $(Z',Z_n) \in B^{n}(0,r_p)$. 
\end{lemma} 

\begin{proof}
The first inequality is immediate. To prove the second inequality, first we choose a
neighbourhood $V$ of $p, \,V \Subset \mathcal{N}_p \cap U^{*}$ such that 
$ \text{diam}(V) < \text{dist} (\overline{V}, \C^n \setminus \mathcal{N}_p) $. Then
$$
  \Sdist(z) \,=\, \delta_{{\mathcal{N}_p} \, \cap \, D} (z) \quad \forall z \in V\cap D.
$$
The choice of $V$ ensures that if $z \in V$, there exists 
$ w_z \in \bdy D \cap \mathcal{N}_p$ such that $\Sdist (z) \,=\, \|z - w_z  \|  $.
\smallskip

Now, clearly the function $Y$ is Lipschitz with some Lipschitz constant $C >1$ and it
vanishes on $ {\sf U}^p (\bdy D \cap \mathcal{N}_p)$. Thus, denoting 
$Z\,=\, {\sf U}^{p}(z)$, we get
\begin{align*}
  Y(Z)\,=\, \| Y({\sf U}^{p}(z))-Y({\sf U}^{p}(w_z))\| \leq C 
  \| {\sf U}^{p}(z) - {\sf U}^{p} (w_z) \| \,=\, C \|z-w_z\| \\ =\, C \Sdist(z) 
  \quad \forall z \in V \cap D.
\end{align*}
Hence the result.
\end{proof}

The next result is the first step on the path\,---\,via the distance decreasing property
for the Kobayashi metric\,---\,to get an integrable bound on the norm of the total
derivative of the map $F$ in the proof of Theorem~\ref{th:non_smooth_bdy_extn}. This will
allow us to use a type of ``Hardy--Littlewood trick'' to establish local
continuous extension of proper holomorphic maps. The latter idea, which we will
adapt to our low-regularity setting, is
inspired by the proof of {\cite[Lemma~8]{diederich-fornaess:1979}} (where the relevant 
estimates are absent, but the idea for obtaining them is hinted at).

\begin{result}[paraphrasing {\cite[Proposition~1.4]{mercer:1993}}]  \label{r:loc_Hopf}
Let $D \subset \Cn$ be a domain and $W$ be an open set such that $\bdy D \cap W$ is a
non-empty $\bdy D$-open set. Suppose $D$ satisfies a uniform interior cone condition in
$W$. Let $\varphi: D \rightarrow [-\infty,0) $ be a plurisubharmonic function. Then, given
any open set $\opu \subset D$ with $\opu \Subset W$, there exist constants
$C_\opu >0$ and $\alpha_{\opu}>1$ such that
\begin{equation}\label{eqn:H_ineq}
  \varphi(w) \leq - C_{\opu} (\Sdist(w))^{\alpha_\opu} \quad \forall w \in \opu.
\end{equation}
\end{result} 

\begin{remark}
We remark here that the above is a local version of the Hopf Lemma due to Mercer given by 
\cite[Proposition~1.4]{mercer:1993}. The proof of the local version is routine and we shall skip it. 
(Indeed, a lemma of Hopf-type is typically a local statement; Mercer's result was
formulated as a global statement presumably because, then, its proof is similar to that
of the global statement given by \cite[Proposition~12.2]{fornaess-stensones:1987}.)
\end{remark}

The next lemma is an application of the above version of a Hopf-type lemma\,---\,suited for the
geometry of $W \cap \Omega$\,---\,applied to the ``pushforward'' of the plurisubharmonic
function $\rho$ on $D$. To be precise, (using classical results for proper holomorphic
maps) the function 
$$ 
  \tau:\Omega \rightarrow \R, \ \tau(w)\,:=\, \text{max} \, \{\rho(z): z \in F^{-1} \{w\} \}, \ w \in \Omega,
$$ 
is a continuous, negative, plurisubharmonic function on $\Omega$. Having this we now prove

\begin{lemma} \label{l:target_Hopf}
Let $U^{*}$ be the open set and $F:D \rightarrow \Omega$ be the map occurring in
Theorem~\ref{th:non_smooth_bdy_extn}. Let $\tau$ be as defined above. Then, there exist
constants (which depend on $U^{*}$) $\alpha _{*} >1$ and $C_{*}>0 $ such that 
$$
  \tau(w)\leq - C_{*}{(\Tdist(w))}^{ {\alpha}_{*}} \quad \forall w \in F (U^{*} \cap D). 
$$
\end{lemma}
\begin{proof}
Let $W^*$ be an open set such that $F(U^{*} \cap D) \Subset W^* \Subset W $.
Since $\Omega$ satisfies a uniform interior cone condition in $W$, 
%
by Result~\ref{r:loc_Hopf}, there exist $C_*\,:=\,C_{W^* \cap \Omega}>0$ and $\alpha_*\,:=\,
{\alpha}_{W^* \cap \Omega}>1 $ such that
$$ 
  \tau (w) \leq - C_* {(\Tdist(w))}^{\alpha_*} \quad \forall w \in W^* \cap \Omega.
$$
This proves the lemma.
\end{proof}

\section{Some geometric lemmas}     \label{sec:geom_lemmas}
In this section we wish to obtain an estimate on $\DKobd$, where $D$ is as in Theorem~\ref{th:Dini_C2-cvx},
similar to the estimate provided by Forstneri\v{c}--Rosay in \cite[Proposition~2.5]{forstneric-rosay:1987}.
Their estimate is obtained by embedding a specific
simply-connected bounded planar domain $\mathcal{D}$ into the domain they considered in
\cite[Proposition~2.5]{forstneric-rosay:1987}. 
The boundary of the domain that they considered is of class $\mathcal{C}^{1,\epsilon}$ near a 
given boundary point. However, in our situation the boundary near $p \in \bdy D$ has lower (namely,
$\Dini$) regularity, so we need to modify their construction. For this reason, we
introduce a class of domains $\plDom$, $\beta,\eps>0$ (see the definition below)\,---\,analogous
to $\mathcal{D}$ in \cite[Proposition~2.5]{forstneric-rosay:1987}\,---\,that can be embedded into
$D$ for a suitable choice of $(\beta, \eps)$. These domains appeared in the work of Maitra
\cite{maitra:2020}. In fact, we shall adapt some of the arguments in
\cite[Proposition~4.2]{maitra:2020} to the present case.
\smallskip

In this discussion we need some definitions. Let $\omega, \ \varphi_p $, and $r_p$ be as introduced by
Definition~\ref{defn:Lipz}-$(2)$. Here $n \geq 2$. Let us define $\omega_p: [0, 2\sqrt{2}r_p) \rightarrow [0, \infty)$ as
\begin{equation}    \label{eqn:mod_cty}
  \omega_p(r) \,:=\, \sup \big\{ \|\nabla \varphi_p(x)- \nabla \varphi_p(y)\|: x,\ y \in B^{n-1}(0, r_p) 
  \times (-r_p,r_p), \  \|x-y\|\leq r \big\}.
\end{equation}
One can check that $\ome_p$ satisfies the following properties:
\begin{itemize}
  \item $\omega_p$ is monotone increasing.
  \smallskip
  \item For any $r \in [0, 2\sqrt{2}r_p), \ \ome_p(r)\leq \ome(r)$. In particular,  $\ome_p$ satisfies a Dini condition.
  \smallskip
  \item $\ome_p$ is sub-additive, i.e., given $\sigma, \ \tau \geq 0, \ \sigma+\tau < 2\sqrt{2}
        r_p, \ \ome_p(\sigma+\tau) \leq \ome_p(\sigma)+\ome_p(\tau)$.
\end{itemize}
Now we define $h:(-2\sqrt{2}r_p,2\sqrt{2}r_p) \rightarrow [0, \infty)$ as 
\begin{align}   \label{eqn:integral_mod_cty}
  h(t)\,:=\,\begin{cases}
  			  \,\displaystyle{\int\limits_{0}^{t}}\omega_p(r)dr,   &\text{if $t \geq 0$}, \\
  			  \,\displaystyle{\int\limits_{t}^{0}} \omega_p(-r) dr, &\text{if $t <0$}.
  			\end{cases}
\end{align}
It is easy to see that $h$ is strictly increasing on $[0,2\sqrt{2}r_p)$, strictly decreasing on 
$(-2 \sqrt{2}r_p,0]$, and $h(0)=h'(0)=0$. Now, given $\beta,\ \eps>0$, define the domain
$$
  \plDom\,:=\, \{\zeta=s+it \in \C: |t|<\eps,\ \beta h(t)<s<\eps \}. 
$$
With these definitions we now prove
\begin{proposition}  \label{prpn:pl_dom_embd}
Let $D,\ U$ be as in Theorem~\ref{th:Dini_C2-cvx}. Let $p \in \bdy D \cap U$ and $\mathcal{N}_p$ 
be as given by Definition~\ref{defn:Lipz}-$(2)$. For $\xi \in \bdy D \cap U$, let $\Psi_{\xi}:\C
\rightarrow\Cn$ denote the $\C$-affine map $\zeta \mapsto \xi + \zeta \nrml$ (where $\nrml$ denotes
the unit inward normal vector at $\xi$). Then, for any neighbourhood $V$ of $p$, $V \Subset
\mathcal{N}_p$, there exist constants (that depend on $V$) $\beta,\ \eps >0$ such that 
$ \Psi_{\xi}(\plDom) \subset U \cap D$ for all $\xi \in \bdy D \cap V$.
\end{proposition}

\begin{proof}
Let ${\sf U}^p, \ \varphi_p $ be as in Definition~\ref{defn:Lipz}-$(2)$. Let us
denote $\widetilde{S}\,=\,{\sf U}^p(S)$ for any subset $S \subset \Cn$. We shall
indicate that we are working in the coordinate system given by ${\sf U}^p$ by using $Z_1, \dots,Z_n$,
where $ \left(Z_1,\dots, Z_n \right) =Z={\sf U}^p(z)$. Clearly, $\rho(Z)\,:=\,
\varphi_p(Z', \re Z_n)- \im Z_n, \ Z\,=\,(Z',Z_n) \in \widetilde{\mathcal{N}_p}$, is a local
defining function for $\widetilde{D}$ near $0\in \bdy \widetilde{D}$. Since $\omega_p$ as defined 
in \eqref{eqn:mod_cty} is increasing, we~have
\begin{equation} \label{eqn:mod_cty_estm}
  \| \nabla \rho(Z)- \nabla \rho(W)\| \leq \omega_p \big(\|Z-W\|\big) \quad \forall Z,\ W \in \widetilde{\mathcal{N}_p}.
\end{equation}

Since $V \Subset \mathcal{N}_p$, $m \equiv m_V \,:=\,\inf \big\{\|D\rho\big(\widetilde{\xi}\, \big)\|:
\widetilde{\xi} \in \bdy \widetilde{D} \cap \widetilde{V} \big\}>0$. Now choose $\beta \equiv \beta_V >1$
such that $1/ \beta \leq m\big/ 4\sqrt{2}$. The need for such a choice for $\beta$
will become evident later. Take $r_V>0$ with $r_V < r_p$ satisfying $ \big( \bigcup_
{\xi \in \bdy D \cap V} B^n(\xi, r_V) \big) \cap D \subset \mathcal{N}_p \cap D$. Also, recall that
the function $h$ as defined in \eqref{eqn:integral_mod_cty} satisfies $h'(0)=0$. So, we can choose 
a sufficiently small $\eps \equiv \eps_V >0$ such that 
\begin{align}   \label{eqn:h'_vanish}
  \sqrt{2}\eps < r_V \quad \text{and} \quad x/ h^{-1}(x) < 1/\beta \ \; \forall x \in (0,\eps).
\end{align}

It is clear that with the choice of $\beta, \ \eps$ above, $\plDom \subset D(0,r_V)$.
It suffices to show that ${\sf U}^p \big( \Psi_{\xi}(\plDom) \big) \subset
\widetilde{\mathcal{N}_p} \cap \widetilde{D}$ for all $\xi \in \bdy D \cap V$. Fix $\xi \in 
\bdy D \cap V$, $\zeta=s+it \in \plDom$. Now, ${\sf U}^p(\xi+\zeta \nrml)\,=\, \widetilde{\xi}+
\zeta \widetilde{\nrml}$ where we write $\widetilde{\xi}\,:=\, {\sf U}^p(\xi)$, 
$ \widetilde{\nrml}\,:=\, {\sf U}_p(\nrml)$. Then, by the Fundamental Theorem of Calculus,
\begin{align*}
  \rho\big(\widetilde{\xi}+ \zeta \widetilde{\nrml}\big) &= \rho \big(\widetilde{\xi}\, \big) +
  D\rho\big(\widetilde{\xi}\,\big) \big(\zeta \widetilde{\nrml} \big) + \int\limits_{0}^{1}
  \big( D\rho\big(\widetilde{\xi}+x \zeta \widetilde{\nrml}\big) -D\rho 
  \big(\widetilde{\xi} \, \big) \big) \big(\zeta \widetilde{\nrml} \big) dx \notag \\
  &=\,s \big\langle \nabla \rho \big(\widetilde{\xi}\,\big) \mid \widetilde{\nrml} \big\rangle + t\big\langle \nabla
  \rho \big(\widetilde{\xi}\,\big)\mid \mathbb{J\big(\widetilde{\nrml}\big)}\big\rangle +
  \int\limits_{0}^{1} \big\langle \nabla\rho \big( \widetilde{\xi}+ x\zeta \widetilde{\nrml} \big) - 
  \nabla \rho \big( \widetilde{\xi}\, \big) \mid \zeta \widetilde{\nrml} \big\rangle \, dx,
\end{align*}
where $\langle\,\bcdot\mid \bcdot\,\rangle$ denotes the standard inner product on
${\R}^{2n}$ and $\mathbb{J}:{\R}^{2n} \rightarrow {\R}^{2n}$ is the standard almost complex
structure. Now, ${\sf U}_p$ being unitary, $\widetilde{\nrml}\,=\, - \nabla \rho \big(\widetilde{\xi} \,\big)
\big/ \| \nabla \rho \big(\widetilde{\xi} \,\big) \|$. Thus, using \eqref{eqn:mod_cty_estm}, the last
equation gives 
\begin{equation}    \label{eqn:Taylor_simplfd}
  \rho \big(\widetilde{\xi}+ \zeta \widetilde{\nrml} \big) \leq -sm + |\zeta| \int\limits_{0}^{1}
  \omega_p (x |\zeta|)\,dx.
\end{equation}
Since $\zeta= s + it \in \plDom$, in view of \eqref{eqn:h'_vanish}
$$
  \ |\zeta|^2 \,\leq\, s^2 + \big(h^{-1}(s/\beta)\big)^2 \,=\, \big(h^{-1}(s/\beta) \big)^2 \bigg( 
  {\beta}^2 \Big(\frac{(s/\beta)}{h^{-1}(s/\beta)}\Big)^2+~1 ~\bigg) \leq 2 \big(h^{-1}(s/\beta) \big)^2.
$$
In view of this estimate, \eqref{eqn:Taylor_simplfd} implies 
\begin{align*}
  \rho \big( \widetilde{\xi} + \zeta \widetilde{\nrml}\big) &\leq -sm + \sqrt{2}
  h^{-1}(s/\beta) \int\limits_{0}^{1} \omega_p \big(2xh^{-1}(s/\beta)\big)\,dx \\ 
  &\leq -sm + 2\sqrt{2} h^{-1}(s/\beta) \int\limits_{0}^{1} \omega_p \big(xh^{-1}(s/\beta) \big)\,dx \\
  &= -sm + 2\sqrt{2}\!\!\!\int\limits_{0}^{h^{-1}(s/\beta)}\!\!\!\omega_p(u)\,du
  \,=\, -sm + 2\sqrt{2}s/\beta  \,\leq\, -sm + sm/2 \,<\, 0,
\end{align*}
where the second inequality is due to the sub-additivity of $\omega_p$ and the last inequality
 is due to our choice of $\beta$. This shows that $\widetilde{\xi} + \zeta \widetilde{\nrml} \in
 \widetilde{\mathcal{N}_p} \cap \widetilde{D}$ for all $\zeta \in \plDom$. Hence the result.
\end{proof}

In order to present the next two results, we require a definition and some related remarks.

\begin{definition}  \label{defn:Riemann_map_diffeo}
A bounded domain $\mathcal{D} \varsubsetneq \C$ is called a \emph{model domain} if $\mathcal{D}$ is an open 
subset of $\{\zeta \in \C: \re \, \zeta > 0 \}$ that is symmetric about $\R$, whose boundary is a Jordan
curve with $0 \in \bdy  \mathcal{D}$, and such that, if $g$ denotes the unique biholomorphic mapping of 
$\mathcal{D}$ onto $\mathbb{D}$ such that $g (\mathcal{D} \cap \R)\,=\, (-1,1)$ and $g(0)\,=\,1$, then
the limit
\begin{equation}\label{eqn:std_domain}
  \lim\limits_{\mathcal{D} \ni \zeta \rightarrow 0} \frac{g(\zeta) -1 }{\zeta}
\end{equation}
exists and is non-zero.
\end{definition}
\begin{remark}
Since $\mathcal{D}$ in Definition~\ref{defn:Riemann_map_diffeo} is enclosed by a Jordan curve, it follows
from Carath\'{e}od-ory's theorem that any biholomorphic map $G$ of $\mathcal{D}$ onto $\mathbb{D}$ extends to
a homeomorphism from $\overline{\mathcal{D}}$ to $\overline{\mathbb{D}}$. Thus, for such a map,
$G(0)$ makes sense. Now, it is classical that a biholomorphic map $g: \mathcal{D} \xrightarrow{\text{onto}}
\mathbb{D}$ such that $g(\mathcal{D} \cap \R)\,=\,(-1,1)$ and $g(0)\,=\,1$ exists and is unique.
\end{remark}

\begin{remark}  \label{rem:std_dom_exmp}
The domains $\plDom, \ \beta, \ \eps>0$, introduced above are examples of model domains. It is obvious 
from its construction that $\plDom$ has the geometric properties that a model domain must have. As for
the existence of the limit \eqref{eqn:std_domain} for the biholomorphic map
$g_{\scriptscriptstyle{\beta},\eps}: \plDom \rightarrow \mathbb{D}$ that
maps $\plDom \cap \R$ onto $(-1,1)$ and such that $g_{\scriptscriptstyle{\beta},\eps}(0)\,=\,1$ and
the condition on it: the desired
condition is given by \cite[Lemma~4.3]{maitra:2020} by Maitra combined with 
\cite[Theorem~1]{warschawski:1961} by Warschawski.
\end{remark}

The next proposition will lead to our desired estimate on $K_D$ as discussed in the beginning of
this section (for $D$ as in Theorem~\ref{th:Dini_C2-cvx}). Although this proposition is stated
for domains that satisfy very general conditions near a given boundary point, we shall apply it
to domains that satisfy a much more geometrical hypothesis. To be explicit, the domains that are
considered in Theorem~\ref{th:Dini_C2-cvx} will be shown to satisfy the hypothesis of the
following proposition. Here $\nrml$ will denote the unit inward normal vector to $\bdy D$ at $\xi \in 
\bdy D \cap V$.

\begin{proposition} \label{prpn:F-R_estm_analg}
Let $D$ be a domain in $\Cn$ and $p \in \bdy D$. Suppose $D$ admits a pair of balls $B^n(p, 4\varrho)
\,=:\,V$ and $B^n(p, \varrho) \,=:\,V'$ such that $\bdy D \cap V$ is a connected $\mathcal{C}^1$-submanifold
of $V$, and such that
\begin{itemize}
  \item $\|\nrml - \eta_p \| < 1/8 $ for every $\xi \in \bdy D \cap V.$
  \smallskip
  
  \item For each $z \in V' \cap D$, there exists a point $\xi \in \bdy D \cap V$ such that $\Sdist(z)\,=\,
        \|z-\xi\|$.
  \smallskip        
        
  \item There exists a constant $c \in (5/8, 1)$ such that 
        $$
          z+ t \nrml \in D \ \ \text{and} \ \ \Sdist(z+t \nrml)> ct
        $$
        for every $t \in [0, 2 \varrho]$, for every $z \in V' \cap D$, and every $\xi \in \bdy D \cap V$.
\end{itemize}
For $\xi \in \bdy D \cap V$, let $\Psi_{\xi}: \C \rightarrow \Cn$ denote the $\C$-affine map
$\zeta \mapsto \xi + \zeta \nrml$. Assume that there exists a model domain $\mathcal{D} \subset \C$ 
such that $\Psi_{\xi}(\mathcal{D}) \subset D$ for each $\xi \in \bdy D \cap V$. Then, there exists a constant 
$C>0$ such that for each $z_1,\ z_2 \in V' \cap D$,
$$
  K_D (z_1,z_2) \leq \sum\limits_{j=1}^{2} \frac{1}{2} \log \frac{1}{\Sdist(z_j)} - \sum\limits_{j=1}^{2}
  \frac{1}{2} \log \bigg(\frac{1}{\Sdist(z_j)+ \|z_2-z_1\|}\bigg) + C.
$$
\end{proposition} 

The above result is, in essence, \cite[Proposition~2.5]{forstneric-rosay:1987} by Forstneri\v{c}--Rosay with
the conditions that make the proof of the latter work emphasised in the hypothesis of 
Proposition~\ref{prpn:F-R_estm_analg} in place of the local $\mathcal{C}^{1,\epsilon}$ condition of
Forstneri\v{c}--Rosay. Thus, the proof of Proposition~\ref{prpn:F-R_estm_analg} is nearly
verbatim the proof of \cite[Proposition~2.5]{forstneric-rosay:1987} (keeping careful track
of where $U$ and $\widetilde{U}$ are required, which correspond to $V'$ and $V$,
respectively, in our case). Therefore, we shall omit the proof of the above proposition.
\smallskip

The following is the estimate on $K_D$ alluded to at the beginning of this section.

\begin{corollary}   \label{cor:source_dom_upp_bd}
Let $D \subset \Cn$ be a domain, $p \in \bdy D$, and $U$ be a neighbourhood of $p$ such that 
$\bdy D \cap U$ is a $\Dini$ submanifold of $U$. Then, there exists a neighbourhood $U'$ of
$p$ and a constant $C > 0$ such that for each $z_1,\ z_2 \in U' \cap D$,
$$
  K_D (z_1,z_2) \leq \sum\limits_{j=1}^{2} \frac{1}{2} \log \frac{1}{\Sdist(z_j)} - \sum\limits_{j=1}^{2}
  \frac{1}{2} \log \bigg(\frac{1}{\Sdist(z_j)+ \|z_2-z_1\|}\bigg) + C.
$$
\end{corollary}

\begin{proof}
The set $\mathcal{N}_p$ below will be as given by Definition~\ref{defn:Lipz}-$(2)$.
Since $\bdy D \cap U$ is a $\Dini$ submanifold of $U$, we can choose a sufficiently small $\varrho >0$
such that $B^n(p, 4 \varrho) \Subset \mathcal{N}_p$ and such that\,---\,writing $V\,:=\, B^n(p, 4 \varrho)$,
$V'\,:=\,B^n(p, \varrho)$\,---\,all but the last condition stated in
Proposition~\ref{prpn:F-R_estm_analg} are satisfied. Now, by Proposition~\ref{prpn:pl_dom_embd}, there
exist constants $\beta, \ \eps>0$ such that $\Psi_{\xi}(\plDom) \subset U \cap D$ for all $\xi \in 
\bdy D \cap V$. Thus, in view of Remark~\ref{rem:std_dom_exmp}, the result follows immediately from 
Proposition~\ref{prpn:F-R_estm_analg} (by taking $U'\,=\,V'$).
\end{proof}

\section{Estimates near a locally log-type convex boundary}  \label{sec:lemmas_Dini_C2-cvx}
This section is devoted to obtaining a useful lower bound for $\TKobd$, where $\Omega$ is as in Theorem~
\ref{th:Dini_C2-cvx}. Such a lower bound is the delicate part of the proof of Theorem~\ref{th:Dini_C2-cvx}.
To this end, we first state a couple of results from the literature.

\begin{result}[Graham, {\cite{graham:1990, graham:1991}}]  \label{r:Graham_cvx}
Let $\Omega$ be a bounded convex domain in $\Cn$. 
Then:
$$
  \frac{\|v\|}{2 \DiscRadi} \leq \TKobm(z;v) \leq \frac{\|v\|}{\DiscRadi} \quad 
  \forall z \in \Omega \ \text{and} \ \forall v \in \Cn \setminus \{0\}.
$$
\end{result}

\noindent{Here, the quantity $\DiscRadi$ is as introduced in Section~\ref{sec:intro}.}

\begin{result}[paraphrasing {\cite[Theorem~3.2]{liu-wang:2021}} by Liu--Wang]  \label{r:L-W_Kobm_loc}
Let $\Omega$ be a bounded domain, $q \in \bdy \Omega$, and suppose there exists a neighbourhood $\opn$
of $q$ such that $\opn \cap \Omega$ is log-type convex. Then, there exists a neighbourhood $W$ of $q$,
$W \Subset \opn$, and a constant $C_q>1$ such that 
$$
  k_{\opn \cap \Omega}(w;v) \leq C_q \TKobm(w;v) \quad \forall w \in W \cap \Omega \ \text{and} \ 
  \forall v \in \Cn.
$$
\end{result}

\begin{remark}  \label{rem:L-W_loc}
{\cite[Theorem~3.2]{liu-wang:2021}} has a seemingly more technical statement than the above
paraphrasing. We get the above\,---\,in the language of {{\cite[Theorem~3.2]{liu-wang:2021}}}\,---\,by focusing 
attention to $\xi\,=\, q \in \bdy \Omega \cap \opn$. Then, $W$ is the ball $B^n(q, \eps)$, where $\eps > 0$
is as given by the latter theorem. The constant $C_q >1$, then, is 
$$
  \exp \left(c {\bigg( \log \frac{1}{\sup_{x \in W \cap \Omega}\Tdist(x)}\bigg)}^{-(1+ \nu)} \right),
$$
where $c, \nu >0$ are as given by {{\cite[Theorem~3.2]{liu-wang:2021}}}.
\end{remark}

\begin{lemma}   \label{l:LTC_Kobm_lower_bd}
Let $\Omega, \ q \in \bdy \Omega$, and $\opn \varsubsetneq_{\rm open} \Cn$ be as in the statement of 
Result~\ref{r:L-W_Kobm_loc}. Then, there exists a neighbourhood $W$ of $q$, $W \Subset \opn$, and 
constants $c, \ \nu >0$ such that
$$
  \TKobm(w;v) \geq c \|v\| {\bigg( \log \frac{1}{\Tdist(w)}\bigg)}^{1+ \nu} \quad 
  \forall w \in W \cap \Omega \ \text{and} \ \forall v \in \Cn.
$$
\end{lemma}

\begin{proof}
Since $\opn \cap \Omega$ is log-type convex, there are constants $C, \nu >0$ such that
$$
  \delta_{\opn \cap \Omega}(w;v) \leq C {\big|\log \delta_{\opn \cap \Omega}(w)\big|}
  ^{-(1+ \nu)} \quad \forall w \in \opn \cap \Omega \ \text{and} \ \forall  v \in \Cn.
$$

We now choose a sufficiently small neighbourhood $W$ of $q$ with $W \Subset \opn $ such that the
conclusion of Result~\ref{r:L-W_Kobm_loc} holds and such that $\delta_{\opn \cap \Omega}(w) \,=\,
\Tdist(w) < 1$ for all $w \in W \cap \Omega $. Then, applying Result~\ref{r:Graham_cvx} to
$k_{\opn \cap \Omega}$ and using the inequality above, the result follows.
\end{proof}

The next proposition relies on the convex domain $\opn \cap \Omega$, where $\opn$ and $\Omega$ 
are as introduced above, having a geometry that is favourable for the estimates that we require.
At this stage, we need to introduce a new notion; we wish to have a brief discussion about
Goldilocks domains, which were introduced by Bharali--Zimmer \cite[Definition~1.1]{bharali-zimmer:2017}. 
In the proposition below, a key step is to show that $\opn \cap \Omega$ is a Goldilocks domain. This
fact and the convexity of $\opn \cap \Omega$ together will allow us to prove this result. First, we need
a couple of definitions:
\smallskip

Let $\Omega \subset \Cn$ be a bounded domain and $r>0$. Define
$$
  M_{\Omega}(r) \,:=\, \sup \bigg\{ \frac{1}{k_{ \Omega}(w;v)} : \delta_{\Omega}(w) \leq r,\ \|v\|\,=\,1 \bigg\}.
$$
We say that $\Omega$ is a \emph{Goldilocks domain} if
\begin{enumerate}
    \item for some (hence any) $\eps >0$ we have
          $$
            \int\limits_{0}^{\eps} \frac{1}{r}M_{\Omega}(r)dr < \infty, \ \text{and}
          $$
    \item for each $w_0 \in \Omega$ there exist constants $\alpha, \ \beta>0$ (that depend on $w_0$)
          such that 
          $$
            \TKobd(w,w_0) \leq \alpha + \beta \log \frac{1}{\Tdist(w)} \quad \forall w \in \Omega.
          $$
\end{enumerate}
See \cite[Remark~1.3]{bharali-zimmer:2017} for an explanation of the geometric significance of the
two conditions above.
\begin{proposition} \label{prpn:target_dom_lower_bd_auxl}
Let $\Omega$, $q \in \bdy \Omega$, and $\opn \varsubsetneq_{\rm open} \Cn$ be as in the statement of
Result~\ref{r:L-W_Kobm_loc}. Then, for $\xi \in (\bdy \Omega \cap \opn) \setminus \{q\}$, there exist
constants $\eps, \ K>0$ such that $B^n(q, \eps), \ B^n(\xi, \eps) \subset \opn$ and
\begin{equation}    \label{eqn:target_dom_lower_bd_auxl}
  K_{\opn \cap \Omega}(w_1,w_2) \geq \frac{1}{2} \log \frac{1}{\Tdist(w_1)} + 
  \frac{1}{2} \log \frac{1}{\Tdist(w_2)} - K
\end{equation}
for all $w_1 \in B^n(q, \eps)\cap \Omega$, for all $w_2 \in B^n(\xi, \eps)\cap \Omega$.
\end{proposition}

\begin{proof}
This result will be proved in the following two steps.

\medskip
\noindent{{\textbf{Step 1.}} \emph{Showing that $\opn \cap \Omega$ is a Goldilocks domain}}
\smallskip

\noindent{First, we shall prove that 
$$
  \int\limits_{0}^{\eps}\frac{1}{r}M_{\opn \cap \Omega}(r) dr < \infty
$$
for some (hence for any) $\eps>0$. Fix $0<\eps<1, \ r \in (0,\eps]$, and $v \in \Cn $ with
$\|v\|\,=\,1$. Let $w \in \opn \cap \Omega$ be such that
$\delta_{\opn \cap \Omega}(w) \leq r$. By Result~\ref{r:Graham_cvx} we get 
\begin{align*}
  \frac{1}{k_{\opn \cap \Omega}(w;v)} \leq 2 \delta_{\opn \cap \Omega}(w;v) \leq 2 C 
  {\bigg( \log \frac{1}{\delta_{\opn \cap \Omega}(w)}\bigg)}^{-(1+ \nu)} \! \! \! \leq 
  2C {\bigg(\log \frac{1}{r} \bigg)}^{-(1+ \nu)},
\end{align*}
where $C, \ \nu >0$ are as given by Definition~\ref{defn:L-T-C}. Since $\int_{0}^{\eps} r^{-1}
{\big(\log(1/r) \big)}^{-(1+ \nu)} dr < \infty$, our desired integral is also convergent.}
\smallskip

Now, $\opn \cap \Omega$ being convex, by \cite[Lemma~2.3]{bharali-zimmer:2017}, for each $w_0 \in \opn \cap
\Omega$ there are constants $\alpha, \ \beta>0$ (that depend on $w_0$) such that
$$
  K_{\opn \cap \Omega} (w_0,w) \leq \alpha + \beta \log \frac{1}{\delta_{\opn \cap \Omega}(w)} \quad 
  \forall w \in \opn \cap \Omega.  
$$

The above argument shows that $\opn \cap \Omega$ satisfies both of the conditions in
\cite[Definition~1.1]{bharali-zimmer:2017} for being a Goldilocks domain.

\pagebreak
\medskip
\noindent{{\textbf{Step 2.}} \emph{Proving the estimate \eqref{eqn:target_dom_lower_bd_auxl}}}
\smallskip

\noindent{Since $\opn \cap \Omega$ is convex, it follows from \cite[(2.4)]{bracci-nikolov-thomas:2022} that
\begin{equation}    \label{eqn:cvx_Kobm_lower_bd}
  K_{\opn \cap \Omega}(w,w') \geq \frac{1}{2} \bigg| \log \frac{\delta_{\opn \cap \Omega}(w)}
  {\delta_{\opn \cap \Omega}(w')} \bigg| \quad \forall w, \ w' \in \opn \cap \Omega.
\end{equation}}

Now, fix a point $o \in \opn \cap \Omega$ and let $\xi \in (\bdy \Omega \cap \opn) \setminus
\{q\}$. So, $q$ and $\xi$ are two distinct boundary points of the Goldilocks domain  $\opn \cap \Omega$.
Let $\eps >0$ be such that $B^n(q, \eps), \ B^n(\xi, \eps) \subset \opn$, such that\,---\,denoting
$V_q\,:=\,B^n(q,\eps) \cap {\Omega}, \ V_{\xi}\,:=\,B^n(\xi,\eps)\cap \Omega$\,---\,$\overline{V_q} \cap 
\overline{V_{\xi}} \,=\,\emptyset $, and such that $\delta_{\opn \cap \Omega}(w)\,=\, \delta_{\Omega}(w)< 
\delta_{\opn \cap \Omega}(o)$ for all $w \in V_q \cup V_{\xi}$. Then, it follows from the proof of
\cite[Proposition~6.8]{bharali-zimmer:2017} that there exists a constant $K'>0$ such that
\begin{equation} \label{eqn:Gromov_prod_bd}
  K_{\opn \cap \Omega}(w_1,o)+K_{\opn \cap \Omega}(o, w_2)-K_{\opn \cap \Omega}(w_1,w_2) \leq K' 
  \quad \forall w_1 \in V_q \ \text{and} \ \forall w_2 \in V_{\xi}.
\end{equation}

Therefore, \eqref{eqn:cvx_Kobm_lower_bd} and \eqref{eqn:Gromov_prod_bd} together imply
$$
  K_{\opn \cap \Omega}(w_1,w_2) \geq  \frac{1}{2} \log \frac{1}{\Tdist(w_1)} + 
  \frac{1}{2} \log \frac{1}{\Tdist(w_2)}+\log \delta_{\opn \cap \Omega}(o) - K'
  \quad \forall w_1 \in V_q \ \text{and} \ \forall w_2 \in V_{\xi}.
$$
Thus, choosing $K > \max(0, K'-\log \delta_{\opn \cap \Omega}(o))$, we get the desired conclusion.
\end{proof}

The following will play a key role in the proof of the main result of this section.

\begin{result}[Liu--Wang, {\cite[Theorem~1.4]{liu-wang:2021}}]  \label{r:L-W_Kobd_loc}
Let $\Omega$ be a bounded domain in $\Cn$. Suppose that there exists a connected open set $\opn$ with 
$\bdy \Omega \cap \opn \neq \emptyset$ and $\opn \cap \Omega$ is log-type convex. Then, for any open
set $W$ with $W \Subset \opn$, there exists $K>0$ such that the Kobayashi distance satisfies
$$
  \TKobd(w_1, w_2) \leq K_{\opn \cap \Omega}(w_1,w_2) \leq \TKobd(w_1, w_2) + K \quad 
  \forall w_1, \ w_2 \in W \cap \Omega.
$$
\end{result}

We now have all the ingredients to establish the lower bound for $\TKobd$ that we need.

\begin{proposition}    \label{prpn:target_dom_lower_bd_final}
Let $\Omega$, $q \in \bdy \Omega$, and $\opn \varsubsetneq_{\rm open} \Cn$ be as in the statement of
Result~\ref{r:L-W_Kobm_loc}. Then, for $\xi \in (\bdy \Omega \cap \opn) \setminus \{q\}$, there exist
constants $\eps, \ K>0$ such that $B^n(q, \eps), \ B^n(\xi, \eps) \subset \opn$ and
\begin{equation*}    \label{eqn:target_dom_lower_bd}
  K_{\Omega}(w_1,w_2) \geq \frac{1}{2} \log \frac{1}{\Tdist(w_1)} + 
  \frac{1}{2} \log \frac{1}{\Tdist(w_2)} - K
\end{equation*}
for all $w_1 \in B^n(q, \eps)\cap \Omega$, for all $w_2 \in B^n(\xi, \eps)\cap \Omega$.
\end{proposition}

\begin{proof}
Let us choose $\eps >0$ such that $B^n(q, \eps) \cup B^n(\xi, \eps) \Subset \opn$, such that 
\eqref{eqn:target_dom_lower_bd_auxl} holds for all $w_1 \in B^n(q, \eps) \cap \Omega$ and for all
$w_2 \in B^n(\xi, \eps) \cap \Omega$. By Result~\ref{r:L-W_Kobd_loc}, there exists $K>0$ such that
$$
   K_{\opn \cap \Omega}(w_1,w_2) \leq \TKobd(w_1, w_2) + K \quad  \forall w_1, \ w_2 \in
   \big(B^n(q, \eps) \cup B^n(\xi, \eps) \big) \cap \Omega.  
$$
Thus, using the above estimate, the result follows immediately from Proposition~
\ref{prpn:target_dom_lower_bd_auxl}.
\end{proof}

\section{The proof of Theorem~\ref{th:non_smooth_bdy_extn}}  \label{sec:first_th_proof}
Before we present the proof, we fix some notations that will be used over the course of the proof. 
We will write $\Ftil\,=\, \big(\Ftil_1, \dots, \Ftil_n \big) \,=\, F \circ \left({\sf U}^{p}\right)^{-1}$,
where ${\sf U}^p$ is as introduced in Definition~\ref{defn:Lipz}. For any subset $S \subset \Cn$,
we will write $\widetilde{S}\,=\,{\sf U}^p(S)$, and we will indicate that we are working in the coordinate
system given by ${\sf U}^p$ by using $Z_1, \dots,Z_n$, as in the proof of Proposition~\ref{prpn:pl_dom_embd}.
\smallskip

Further, in the proof below $C$ will denote a positive constant that depends 
\emph{only} on the data in the hypothesis of Theorem~\ref{th:non_smooth_bdy_extn}.
However, the magnitude of $C$ may change from line to line.

\begin{proof}[Proof of Theorem~\ref{th:non_smooth_bdy_extn}]
Let $V$ be the neighbourhood of $p$ as given by Lemma~\ref{l:Lipz_vertical}. Then, since $F(V \cap D) \subset
W \cap \Omega$, for all $z \in V \cap D$ and $v \in \Cn $, we have
$$
  \frac{\|F'(z) v\|}{M \big(\Tdist(F(z)) \big)} \leq \TKobm (F(z);F'(z)v) \leq k_D (z;v)
  \leq  \frac{\|v\|}{\Sdist(z)}.
$$
The last inequality follows by the distance decreasing property for the Kobayashi metric 
under the inclusion map $B^n(z,\Sdist (z)) \hookrightarrow D $. Now, using
Lemma~\ref{l:target_Hopf}, for all $z \in V \cap D$
\begin{equation}\label{eqn:use_of_Hopf}
  -(\Sdist(z))^s \leq \rho(z) \leq \tau(F(z)) \leq - C_* \big( \Tdist(F(z)) \big)^{\alpha_*},
\end{equation}
which gives (since $M$ is increasing), for all $z \in V \cap D$
$$
  M \big( \Tdist(F(z)) \big) \leq M \big( C (\Sdist(z))^{s/{\alpha_*}} \big).
$$
Combining the above inequalities we get
\begin{equation} \label{eqn:F_prime_bd}
  \|F'(z)v\| \leq \frac{\|v\|}{\Sdist(z)}M \big( C(\Sdist(z))^{s/{\alpha_*}} \big) \quad \forall
  z \in V \cap D , \, \forall v \in \Cn.
\end{equation}

We now express the estimate~\eqref{eqn:F_prime_bd} in the new coordinate system $ \Cn \ni z
\mapsto Z= {\sf U}^p(z) $. Let us write $A\, :=\, {\sf U}_{p}^{-1}$, which is a unitary transformation.
Applying chain rule, \eqref{eqn:F_prime_bd} gives
\begin{align} \label{eqn:Ftilde_prime_bd}
  \|\Ftil'(Z)v \|= \|F'(z)(Av) \| \leq \frac{\|Av\|}{\Sdist(z)} M \big( C(\Sdist(z))^{s/{\alpha_*}} \big) 
  &= \frac{\|v\|}{\Sdist(z)} M \big( C(\Sdist(z))^{s/{\alpha_*}} \big) \notag \\ 
  &\leq \frac{C\, \|v\|}{Y(Z)} M \big( C(Y(Z))^{s/{\alpha_*}} \big),
\end{align}
where the last inequality holds by Lemma~\ref{l:Lipz_vertical}. Recall that $M$ satisfies a Dini condition.
It is elementary to see (using change of variables) that for fixed $\kappa,\, m >0$, the 
composite function $t \mapsto M(\kappa\,t^m)$ also satisfies a Dini condition. Hence,
we can rewrite \eqref{eqn:Ftilde_prime_bd} as:
\begin{equation} \label{eqn:Ftilde_prime_bd_final}
  \|\Ftil'(Z)v\| \leq \psi \big(Y(Z) \big) \|v\| \quad \forall Z\in \widetilde{V} \cap
  \widetilde{D}, \, \forall v \in \Cn,       
\end{equation}
where $\psi$ is a non-negative Lebesgue integrable function on ${\rm range}(Y)$.
\smallskip

Now we choose a neighbourhood $\Util$ of $\boldsymbol{0} \in \Cn, \ \Util \Subset \widetilde{V}
$ and $\delta>0$ such that
\begin{align}   \label{eqn:U_tilde}
  \Bigg( \bigcup\limits_{\xi\, \in\, \bdy \widetilde{D} \, \cap\, \Util}
  \overline{B^n (\xi, \delta)} \Bigg) \cap \widetilde{D} \subset \widetilde{V} \cap \widetilde{D}.
\end{align}
Pick $1\leq j \leq n$, fix $\xi=(\xi ', \zeta + i \, \eta) = (\xi', \zeta + i\,
\varphi_p(\xi', \zeta))\in \bdy \widetilde{D} \cap \Util$, and $0<t<t'< \delta$. Then,
\begin{align}  \label{eqn:fund_th}
  \Ftil_j (\xi + t'\ep) - \Ftil_j(\xi + t\ep) = \int\limits_{t}^{t'} i\, \frac{\bdy \Ftil_j}{\bdy Z_n}
  (\xi + x\ep)\, dx \quad \text{where} \ \ep = (0,\dots,0,i) \in \Cn.
\end{align}
By \eqref{eqn:Ftilde_prime_bd_final},
$$
  \Bigg| \int\limits_{t}^{t'} i\, \frac{\bdy \Ftil_j}{\bdy Z_n} (\xi + x\ep)\, dx \Bigg| \leq
  \int\limits_{t}^{t'} \psi \big(Y(\xi+x \ep) \big)\, dx = \int\limits_{t}^{t'}  \psi (x) \, dx.
$$
Hence, by integrability of $\psi$ and the fact that $0 \in$ domain$(\psi)$, the limit 
$$
  \Ftil^{\bullet}_{j}(\xi) \,:=\, \Ftil_j(\xi + t' \ep) - \lim \limits_{t \rightarrow 0^+}
  \int\limits_{t}^{t'} i\, \frac{\bdy \Ftil_j}{\bdy Z_n} (\xi + x\ep)\,dx 
$$
exists. (Note that, although the right hand side of the above expression involves the
parameter $t'$, in view of \eqref{eqn:fund_th}, it does not depend on $t'$.) Also, we get the
estimate (which is uniform in $\xi$)
\begin{align} \label{eqn:bd_int}
  \big| \Ftil_{j}^{\bullet}(\xi) - \Ftil_{j}(\xi + t' \ep) \big| \leq \int\limits_{0}^{t'}
  \psi(x)\,dx \quad \forall \xi \in \bdy \widetilde{D} \cap \Util \ \text{and} \ \forall t' \in (0, \delta).
\end{align}

We can now define $\Fhat\,=\, \big(\Fhat_1, \dots, \Fhat_n \big):\widetilde{D} \cup
\big(\bdy\widetilde{D} \cap \Util \big) \longrightarrow \overline{\Omega}$
\begin{equation}  \label{eqn:ext_map}
  \Fhat(Z)\,:=\,\begin{cases}
                  \Ftil^{\bullet}(Z)\,=\, \big( \Ftil_{1}^{\bullet}(Z), \dots, \Ftil_{n}^{\bullet}(Z) \big),
                  &\text{if $ Z \in \bdy \widetilde{D} \cap \Util$}, \\
                  {} & {} \\
                  \Ftil(Z), &\text{otherwise}.
  			     \end{cases}
\end{equation}
Our goal is to show that $\Fhat$ is continuous. It is enough to show its continuity on $\bdy
\widetilde{D} \cap \Util $. We will adapt a Hardy--Littlewood trick to complete the proof.
\smallskip

Let $\eps>0$. As $\psi$ is integrable near $0$, there exists $r(\eps) \in (0, \delta)$ such that \eqref{eqn:bd_int}
gives
\begin{align}  \label{eqn:eps_bd_int}
   \big| \Ftil_{j}^{\bullet}(\xi) - \Ftil_{j}(\xi + r(\eps) \ep) \big| \leq \int\limits_{0}^{r(\eps)}
   \psi(x)\,dx < \eps/3 \quad \ \text{for every} \ \xi \in \bdy \widetilde{D} \cap \Util.
\end{align}
Define $S(\eps)\,:=\, \left\{ \xi + r(\eps) \ep: \xi  \in \bdy \widetilde{D} \cap \Util \right\}$.
By the choice of $\Util$ as in \eqref{eqn:U_tilde}, we have $S(\eps) \Subset 
\widetilde{V} \cap \widetilde{D} \subset \widetilde{U} \cap \widetilde{D}$. Hence $\Ftil_j
\big|_{\scriptscriptstyle{S(\eps)}}$ is uniformly continuous. So, there exists $\sigma \equiv 
\sigma(\eps) >0$ such that
\begin{align}   \label{eqn:unif_cont}
  \big| \Ftil_j(Z) - \Ftil_j(W) \big| < \eps/3 \quad \text{whenever} \ Z,\,W \in S(\eps) \ 
  \text{with} \ \| Z-W \| < \sigma.
\end{align}
Now, let $\xi_1,\, \xi_2 \in \bdy \widetilde{D} \cap \Util$ such that $\|\xi_1 - \xi_2 \| < \sigma$. 
Let $W_1,\, W_2 \in S(\eps)$ such that $W_k \,:=\, \xi_k + r(\eps) \ep, \ k=1, \, 2$. Since $\|W_1-W_2 \|
= \| \xi_1-\xi_2 \| < \sigma$, \eqref{eqn:eps_bd_int} and \eqref{eqn:unif_cont} together imply
$$
  \big| \Ftil_{j}^{\bullet}(\xi_1) - \Ftil_{j}^{\bullet}(\xi_2) \big| \leq \sum_{k=1}^{2}
  \big| \Ftil_{j}^{\bullet}(\xi_k) - \Ftil_{j}(W_k) \big| + 
  \big| \Ftil_j (W_1) - \Ftil_j (W_2) \big| < \eps.
$$
This shows that $\Fhat \big|_{\scriptscriptstyle{\bdy \widetilde{D} \,\cap\, \Util}}$ is continuous.
\smallskip

Let us now fix $ \xi=(\xi', \zeta + i \, \eta)= (\xi', \zeta+ i \, \varphi_p(\xi',\zeta)) \in 
\bdy \widetilde{D} \cap \Util $. It suffices to show that for any sequence $\{Z_\nu \} \subset 
\Big( \Util \cap \overline{\widetilde{D}} \Big) \setminus \{\xi \}$, converging to $\xi$, we have 
$\Fhat_j(Z_\nu) \rightarrow \Ftil_{j}^{\bullet}(\xi)$ as $\nu \rightarrow \infty$.
Define an auxiliary sequence $\{\widetilde{Z_\nu} \}$ as
$$
  \widetilde{Z_\nu}\,:=\,\begin{cases}
  					Z_\nu, &\text{if $Z_\nu \in \bdy \widetilde{D}$}, \\
  					Z_\nu, &\text{if $Z_\nu= \big(\xi', \zeta + i\, (x+\varphi_p(\xi',\zeta)) \big)$ for some $x>0$}, \\
  					\pi(Z_\nu), &\text{otherwise},
  					\end{cases}
$$
where $\pi (Z)\,:=\, (Z', X_n + i \, \varphi_p(Z',X_n))$ and $Z=(Z',X_n+ i\, Y_n) \in
\Util \cap \widetilde{D} $\,---\,i.e., the projection along $\R\ep$ into the boundary of $\widetilde{D}$.
Since $\pi$ is continuous, $\widetilde{Z_\nu} \rightarrow \xi$ as $\nu \rightarrow \infty$. Now, using
the estimates \eqref{eqn:bd_int} and \eqref{eqn:eps_bd_int} together with the fact that $\Fhat
\big|_{\scriptscriptstyle{\bdy \widetilde{D}\,\cap\,\Util}}$ is continuous, we get
$$
  \lim_{\nu \rightarrow \infty} \Big( \Fhat_j (\widetilde{Z_\nu})- \Fhat_j(Z_\nu) \Big)=0 \ \ \text{and}
  \ \ \lim_{\nu \rightarrow \infty} \Fhat_j(\widetilde{Z_\nu})= \Ftil_{j}^{\bullet}(\xi) = \Fhat_j(\xi).
$$

This proves that $\Fhat$ as defined in \eqref{eqn:ext_map} is continuous. Since ${\sf U}^p$ is an
automorphism of $ \Cn $, this completes the proof.
\end{proof}

Before we end this section, we must elaborate upon a point that was deferred in
Section~\ref{sec:intro}.

\begin{remark}\label{rem:Kobm}
The curious reader might ask whether, since $\Omega$ in Theorem~\ref{th:non_smooth_bdy_extn}
is assumed to satisfy a uniform interior cone condition in $W$, one also requires the
condition \eqref{eqn:Kobm_bound}. The question might arise as the assumption of a uniform
interior cone condition in $W$ suggests the existence of a local plurisuharmonic barrier. Then,
by Result~\ref{r:Kob_low_bd_Sibony}, one may hope to deduce \eqref{eqn:Kobm_bound}. In fact,
the last two ingredients summarise the approach to the proof in \cite{sukhov:1994} (also see
\cite{berteloot-sukhov:1997}). However, such an approach actually requires a local plurisuharmonic
barrier with \emph{``good'' estimates}\,---\,obtaining such estimates is quite hard, and a
uniform interior cone condition is not sufficient for such an estimate. Now, note that the
hypothesis of Theorem~\ref{th:non_smooth_bdy_extn} imposes no regularity condition on
$\bdy{\Omega}$. In its absence, our uniform interior cone condition in $W$ functions as a very
mild (local) boundary-regularity condition on $\Omega$. The latter enables us to use a version of
the Hopf Lemma to obtain the inequality \eqref{eqn:use_of_Hopf}. In contrast, in
\cite{sukhov:1994}\,---\,and in most of the articles on the present theme cited above\,---\,the
relevant patch of the the boundary of the target domain is required to be $\smoo^2$-smooth (in
order to deduce the analogue of \eqref{eqn:use_of_Hopf}). That brings us to 
the condition \eqref{eqn:Kobm_bound}: a careful perusal of \cite{sukhov:1994, berteloot-sukhov:1997}
will reveal that \eqref{eqn:Kobm_bound} is \textbf{less restrictive} than the plurisubharmonic-barrier
condition. Our function $M$ is just a way to express quantitatively the requirement that
$k_{\Omega}(w; \bcdot)$ must grow as $w$ approaches $\bdy \Omega$ via $W \cap \Omega$ but may do so
\textbf{relatively slowly}. Classically, \eqref{eqn:Kobm_bound} is true if $\bdy{\OM} \cap W$
is strongly pseudoconvex (with $M(r) = \sqrt{r}$); it is also true if $\bdy{\OM} \cap W$ is real
analytic and all points in $\bdy{\OM} \cap W$ are of finite type (in which case $M$ is some fractional
power); see \cite{diederich-fornaess:1979}.
However, recall that we want Theorem~\ref{th:non_smooth_bdy_extn} to be able to address the case
when $\bdy{\OM} \cap W$ contains infinite-type points. This is what $M$ (which is more general than
a power) enables; see \cite{bharali:2016, liu-wang:2021}.    
\end{remark} 

\section{The proof of Theorem~\ref{th:Dini_C2-cvx}} \label{sec:second_th_proof}
We will present the proof in two steps. The method that we will use in the first step is 
the method of the proof of \cite[Theorem~1.1]{forstneric-rosay:1987}.
In this step, we will show that $F$ extends as a continuous map on $D \cup \{p\}$. In the
second step, we will apply Theorem~\ref{th:non_smooth_bdy_extn} to get the desired conclusion.

\medskip
\noindent{{\textbf{Step 1.}} \emph{Proving that $F$ extends continuously to $p \in \bdy D$}}
\smallskip

\noindent{First, we note that the cluster set $\cluster$ is connected. This is so because there is a
basis of neighbourhoods $\{{\mathcal{N}}_{\nu}\}$ of $p$ such $\mathcal{N}_{\nu} \cap D$ is
connected for each $\nu$ (see \cite[Chapter~1, Section~1]{collingwood-lohwater:1966} for more
details). Also, $F$ being proper, $\cluster \subseteq \bdy \Omega$.}
\smallskip

We shall establish our goal by contradiction. Assume that $F$ does not extend continuously to 
$p \in \bdy D$. Then, $\cluster$ is not a singleton. Since $\cluster$ is connected,
we can find a point $\xi \in (\bdy \Omega \cap \opn) \cap \cluster$ with $\xi \neq q$. Consider 
a pair of sequences $\{\domseq{1}\}$ and $\{\domseq{2}\}$ in $U \cap D$ such that, writing
$\targseq{j}\,:=\,F\big(\domseq{j}\big)$ for $j\,=\,1, \ 2$, we have
$$
  \domseq{j} \to p \ \text{for} \  j\,=\, 1,\ 2, \ \text{and} \ 
  \targseq{1} \to q, \ \targseq{2} \to \xi \ \text{as $\nu \to \infty$.}
$$
So, we can find a non-negative integer $N>1$ such that for all $\nu \geq N$
\begin{align}
   K_D \big(\domseq{1}, \domseq{2}\big) &\leq \frac{1}{2} \log \frac{1}{\Sdist\big(\domseq{1}\big)}+
   \frac{1}{2} \log \frac{1}{\Sdist\big(\domseq{2}\big)} - l(\nu) + C, \label{eqn:upp_bd_seq}\\
  K_{\Omega}\big(\targseq{1},\targseq{2}\big) &\geq \frac{1}{2} \log
  \frac{1}{\Tdist\big(\targseq{1}\big)} + \frac{1}{2}
  \log\frac{1}{\Tdist\big(\targseq{2}\big)} - K     \label{eqn:low_bd_seq}, 
\end{align}
where $C>0$ is given by Corollary~\ref{cor:source_dom_upp_bd}, $K>0$ is given by Proposition~
\ref{prpn:target_dom_lower_bd_final}, and 
\begin{equation}    \label{eqn:l_nu_unbdd}
  l(\nu)\,:=\, \sum\limits_{j=1}^{2} \frac{1}{2} \log \Bigg(\frac{1}{\Sdist\big(\domseq{j}\big)+
  \big\|\domseq{1}-\domseq{2}\big\|}\Bigg) \to \infty.
\end{equation}
Recall that, the function $\tau:\Omega \rightarrow \R, \ w \mapsto \text{max} \, 
\{\rho(z): z \in F^{-1} \{w\} \},$ is continuous, negative, and plurisubharmonic
on $\Omega$ (see the discussion prior to Lemma~\ref{l:target_Hopf}). Since
$\bdy \Omega$ is of class $\mathcal{C}^2$ near $q$, we appeal to the classical Hopf Lemma
for plurisubharmonic functions which, essentially, is the inequality \eqref{eqn:H_ineq}
wherein\,---\,due to (local) $\mathcal{C}^2$-regularity\,---\,the
exponent $\alpha_{\opu}\,=\,1$ independent of the open set $\opu$. This lemma gives us, raising
the value of the $N$ above if needed, a constant
$C_0>0$ such that for $j\,=1,\ 2$, and for all $\nu \geq N$,
\begin{equation*}   
  \tau\big(\targseq{j} \big) \leq - C_0 \Tdist\big(\targseq{j}\big).
\end{equation*}
This gives, for $j\,=\,1, \ 2 $, and for all $\nu \geq N$,
\begin{equation} \label{eqn:compare_distances}
  -\Sdist\big( \domseq{j} \big) \leq \rho \big( \domseq{j} \big) \leq \tau 
  \big( \targseq{j} \big) \leq -C_0 \Tdist\big( \targseq{j} \big).
\end{equation} 
Then, \eqref{eqn:upp_bd_seq},
\eqref{eqn:low_bd_seq}, and \eqref{eqn:compare_distances} together imply
\begin{align*}
  K_D \big(\domseq{1}, \domseq{2}\big) &\leq \frac{1}{2} \log \frac{1}{\Tdist\big(\targseq{1}\big)}+
  \frac{1}{2} \log \frac{1}{\Tdist\big(\targseq{2}\big)} - l(\nu) + C - \log C_0 \\
  &\leq \TKobd\big(\targseq{1},\targseq{2}\big) + (K+C-\log C_0) - l(\nu) \\
  & \leq K_D \big(\domseq{1}, \domseq{2}\big) + (K+C-\log C_0) - l(\nu) \quad \forall \nu \geq N.
\end{align*}
The above estimate implies that $\{l(\nu): \nu \in \Z_{+}\}$ is bounded, which contradicts
\eqref{eqn:l_nu_unbdd}.
\smallskip

Thus, the map $\Fhat:D \cup \{p\} \rightarrow \overline{\Omega}$ as defined by
\begin{equation}  \label{eqn:extn_at_one_pt}
  \Fhat(z)\,:=\,\begin{cases}
                  F(z), &\text{if $z \in D$}, \\
                  q, &\text{if $z\,=\,p$},
  			     \end{cases}
\end{equation}
is continuous and $\Fhat\big|_{D}\,=\,F$.

\medskip
\noindent{{\textbf{Step 2.}} \emph{Proving that $F$ extends continuously to a $\bdy D$-neighbourhood of $p$}}
\smallskip

\noindent{Let $W$ be the neighbourhood of $q$ with $W \Subset \opn$ and let $C, \ \nu>0$ 
be the constants as given by Lemma~\ref{l:LTC_Kobm_lower_bd}.} Then, 
\begin{equation*}
  \TKobm(w;v) \geq c \|v\| {\bigg( \log \frac{1}{\Tdist(w)}\bigg)}^{1+ \nu} \quad 
  \forall w \in W \cap \Omega \ \text{and} \ \forall v \in \Cn.
\end{equation*}
Note that the function  $r \mapsto {\big(\log(1/r) \big)}^{-(1+ \nu)}, \ 0<r<1$, is
integrable at zero, and its value approaches $0$ as $r \rightarrow 0^{+}$. Thus, it satisfies
the conditions on $M$ featured in Theorem~\ref{th:non_smooth_bdy_extn}.
\smallskip

Now, $\bdy \Omega \cap \opn$ being $\mathcal{C}^2$ smooth, $\Omega$ must satisfy a
uniform interior cone condition in $W$. Also, since the map $\Fhat$ as defined by
\eqref{eqn:extn_at_one_pt} is continuous, we can find an open ball $U^*$ with centre $p$ with
$U^* \Subset U$ such that 
$$
F(U^* \cap D) \subset \Fhat(U^* \cap (D \cup \{p\})) \Subset
W.
$$
Hence, by Theorem~\ref{th:non_smooth_bdy_extn}, $F$ extends continuously to a 
$\bdy D$-neighbourhood of $p$. \hfill \qed
\smallskip

Now that we have seen the argument for Theorem~\ref{th:Dini_C2-cvx}, we can
make the following

\begin{remark} \label{rem:cvx_visi}
We see that the two key properties that the pair $(\Omega, q)$ must possess,
in addition to the issue of the regularity $\bdy \Omega$ near $q$ that was discussed
in Section~\ref{sec:intro}, that make the above proof work are $(i)$~a localization result for
$K_{\Omega}$ akin to Result~\ref{r:L-W_Kobd_loc}, and $(ii)$~whatever that leads to the estimate
\eqref{eqn:Gromov_prod_bd}. Readers familiar with the results in
\cite{bracci-nikolov-thomas:2022} might notice that a weaker condition than
$\opn \cap \Omega$ being log-type convex provides just these two ingredients.
That said:
\begin{itemize}
  \item Replacing log-type convexity of $\opn \cap \Omega$ by the property alluded
  to only ensures the continuous extension of $F$ to $D \cup \{p\}$. It is
  unclear if the continuous extension of $F$ to $D \cup \bdynbd$ can be achieved without
  imposing \emph{some} additional constraint on $k_{\Omega}$.
  \item Our intention in Theorem~\ref{th:Dini_C2-cvx} is to present\,---\,in terms of the
  assumptions made\,---\,a result in the same spirit as Result~\ref{r:F-R_main}. In the
  latter case, a lower bound akin to \eqref{eqn:Kobm_bound} is automatic.
\end{itemize}
It is for these reasons that we highlight Theorem~\ref{th:Dini_C2-cvx} in
Section~\ref{sec:intro}. 
Also, while some conditions are now known that will ensure that
$\opn \cap \Omega$ is convex and has the visibility property (i.e., the property
introduced in \cite{bracci-nikolov-thomas:2022}), the latter property
is slightly non-explicit. Moreover, the latter is not the last word on properties
that will produce a localization result
akin to Result~\ref{r:L-W_Kobd_loc}; this will be discussed in forthcoming work.
We also refer the reader to \cite[Proposition~12]{nikolov-andreev:2017} by
Nikolov--Andreev for a result with assumptions similar to those in
Theorem~\ref{th:Dini_C2-cvx}, but applied globally.
\end{remark}

The discussion in Remark~\ref{rem:cvx_visi} provides enough clues for us to conclude
with one last result. We refer the reader to \cite[Section~2]{bracci-nikolov-thomas:2022}
for definitions.

\begin{theorem}  \label{th:Dini_C2-cvx-visi}
Let $D$ and $\Omega$ be domains in ${\C}^n, \ n\geq 2, \ \Omega$ bounded, and 
let $F: D \rightarrow \Omega$
be a proper holomorphic map. Let $p\in \bdy D$ and $q \in C(F,p)$. Assume that there is a
continuous, negative plurisubharmonic function $\rho$ on $D$ and a neighbourhood $U$ of $p$ 
such that $\bdy D \cap U$ is a $\Dini$ submanifold of $U$, and $\rho(z) 
\geq - \delta_D (z)$ for all $z\in U \cap D$. Suppose there exists a neighbourhood $\opn$
of $q$ such that 
\begin{itemize}
  \item $\bdy \Omega \cap \opn$ is a $\mathcal{C}^2$ submanifold of $\opn$, and
  \smallskip
  \item $\opn \cap \Omega$ is convex and has the visibility property.
\end{itemize}
Then, $F$ extends to a continuous map on $D \cup \{p\}$.
\end{theorem}

Since the proof of the above is very similar to the argument in Step~1 of the
proof in Section~\ref{sec:second_th_proof}, we only mention the differences
in the argument. Firstly, the role of Result~\ref{r:L-W_Kobd_loc} is played
by \cite[Theorem~1.4]{bracci-nikolov-thomas:2022}. Secondly, the visibility
property of $\opn \cap \Omega$ ensures that \eqref{eqn:Gromov_prod_bd} is true.

\section*{Acknowledgements}
\noindent{I am grateful to my thesis advisor, Prof. Gautam Bharali, for suggesting various
ideas and for the invaluable discussions over the course of this project. I am also 
grateful to him for his help with the writing of this paper.
This work is supported by a scholarship from the National Board for Higher Mathematics (NBHM)
(Ref. No. 0203/16(19)/2018-R\&D-II/10706).}

\end{document}